\documentclass[11pt]{article}

\usepackage{amsmath,amssymb,latexsym}
\usepackage{mathrsfs}
\usepackage{verbatim}

\newtheorem{defn}{Definition}[section]
\newtheorem{lem}[defn]{Lemma}
\newtheorem{exmp}[defn]{Example}
\newtheorem{thm}[defn]{Theorem}
\newtheorem{prop}[defn]{Proposition}
\newtheorem{cor}[defn]{Corollary}

\newtheorem{rem}[defn]{Remark}

\newcommand{\ru}{R(U)}
\newcommand{\rus}{R(U_s)}
\newcommand{\uinv}{{U^{-1}}}
\newcommand{\mn}{\mathbb N}
\newcommand{\mr}{\mathbb R}
\newcommand{\mz}{\mathbb Z}
\newcommand{\ma}{\mathcal{A}}
\newcommand{\nulel}{{\bf 0}}
\newcommand{\seq[1]}{\{#1_i\}}
\newcommand{\seqgr[1]}{\{#1_i\}_{i=1}^\infty }
\newcommand{\sech[1]}{{\cap_{s\in\mn_0}{#1}_s}}

\newcommand{\snorm[1]}{{\|\hspace{-.01in}|{#1}\|\hspace{-.015in}|}}
\newcommand{\snormb[1]}{{\left\|\hspace{-.015in}\left|{#1}\right\|\hspace{-.015in}\right|}}
\def\be{\begin{equation}}
\def\ee{\end{equation}}
\def\newin {\kern-0.22em\in\kern-0.15em}
\def\newsubset {\kern-0.2em\subset\kern-0.2em}
\def\vn{\vspace{.1in}\noindent}
\def\v{\vspace{.1in}}
\def\sumii{\sum_{i=1}^{\infty}}
\def\bp{\noindent{\bf Proof. \ }}
\def\ep{\noindent{$\Box$}}
\def\<{\langle}
\def\>{\rangle}
\begin{document}

\title{Series expansions in Fr\'echet spaces and their duals; construction of
Fr\'echet frames}

\author{Stevan Pilipovi\'c\thanks{This research was supported by Ministry of Science of Serbia, Project 144016} \ and Diana T. Stoeva\thanks{This research was supported by DAAD.}}

\maketitle

\begin{abstract}
Frames for Fr\'echet spaces $X_F$ with respect to Fr\'echet
sequence spaces $\Theta_F$ are studied and conditions, implying
series expansions in $X_F$ and $X_F^*$, are determined.
 If $\seqgr[g]$ is a $\Theta_0$-frame for $X_0$, we construct a sequence
$\{X_s\}_{s\in\mn_0}$, $X_s\subset X_{s-1}$, $s\newin \mn$, for
given $\Theta_F$, respectively a sequence
$\{\Theta_s\}_{s\in\mn_0}$, $\Theta_s\subset \Theta_{s-1}$,
$s\newin \mn$, for given $X_F$, so that $\seqgr[g]$ is a
pre-$F$-frame (or $F$-frame) for $X_F$ with respect to $\Theta_F$
under different assumptions given on $X_0$, $\Theta_0$, $\Theta_F$
or $X_0$, $\Theta_0$, $X_F$.
\end{abstract}



\section{Introduction}

In this paper we are interested in the frame expansions of
elements in $\cap_{s=0}^\infty X_s$, where $X_0\supset X_1\supset
X_2\supset\ldots $ are Banach spaces, which corresponds to
$\cap_{s=0}^\infty \Theta_s$, where $\Theta_0\supset
\Theta_1\supset \Theta_2\supset\ldots$ are Banach sequence spaces.
We introduce Fr\'echet frames (cf. \cite{pst}, also) and in this
context we are interested in  the problem of determination of
$\{X_s\}_{s\in\mn_0}$ (resp. $\{\Theta_s\}_{s\in\mn_0}$) for given
frame $\{g_i\}_{i=1}^\infty$ and given $\{\Theta_s\}_{s\in\mn_0}$
(resp. $\{X_s\}_{s\in\mn_0}$). Our approach is different from the
approach of Feichtinger and Gr\" ochening \cite{GFa, FGb} who have
used orbit, respectively, coorbit spaces, for the construction of
a linear space with coefficients in a given sequence space,
respectively, for the construction of a sequence space which
corresponds to a given function space. Our main aim is the
analysis of conditions on $\{g_i\}_{i=1}^\infty$, $\Theta_0$ and
$X_0$ and their subspaces in order to obtain a fine
characterization of  $f$ belonging to a subspace of $X_0$ which
can be a Banach space or a Fr\'echet space, through the properties
of the  corresponding sequence of coefficients
$\{g_i(f)\}_{i=1}^\infty$.

Historically, theory of frames appeared in the paper of Duffin and
Schaeffer \cite{ds} in 1952. Around 1986, Daubechies, Grossmann,
Meyer \cite{DGM} and others reconsidered ideas of Duffin and
Schaeffer, and started to develop the wavelet and frame theory.
Banach frames were  introduced by Gr\"{o}chenig \cite{G} and
subsequently many mathematicians have contributed to this theory,
see for example \cite{AST,ASTproc,Cas,CCS,CH} and references
therein. In the last years, Banach frames for families of Banach
spaces were studied in several papers. Using a tight wavelet frame
for $L^2(\mr^d)$, Borup, Gribonval and Nielsen \cite{BN} obtained
a $\Theta_p$-frame for $L^p(\mr^d),$ for an  appropriate sequence
space $\Theta_p$,
 which gives series expansions in $L^p(\mr^d)$, $1<p<\infty$.
Gr\"{o}chenig and his collaborators \cite{cogr2,forgr,Gr3},
considered localized frames and through the orbit and coorbit
spaces obtained Banach frames for a family of Banach spaces.
Considering $p$-frames for shift-invariant subspaces of $L^p$,
Aldroubi, Sun and Tung \cite{AST} proved that when a sequence of
translations of a finite set of appropriate functions
$\phi_1,...,\phi_r$ forms an $\ell^{p_0}$-frame for the
shift-invariant space $V_{p_0}(\phi)\subset L^{p_0}$, for some
$p_0>1$, then this sequence is also an $\ell^p$-frame for
$V_p(\phi)$ for all values of $p>1$.

In the present paper we consider the following problems:
 \begin{itemize}
 \item[ 1.] Determine conditions, which imply series expansions in projective and inductive limits of Banach spaces.
\end{itemize}
Let $\{(X_s,||\cdot||_s)\}_{s\in\mn_0}$ be a sequence of Banach
spaces such that $X_s\subset X_{s-1}\subset X_0$, $s\newin\mn$, in
the topological sense also, and $X_F=\cap_{s\in\mn_0} X_s$ be
dense in every space $X_s$. Let
$\{(\Theta_s,\snorm[\cdot]_s)\}_{s\in\mn_0}$ be a family of Banach
sequence spaces with similar properties. We determine conditions
on a sequence $\seqgr[g]$, $g_i\newin X_F^*$, which imply the
existence of $\seqgr[f]$, $f_i\newin X_F$, such that every
$f\newin X_F$ and every $g\newin X_F^*$ can be written as
$f=\sum_{i=1}^{\infty} g_i(f) f_i$ and $g=\sum_{i=1}^{\infty}
g(f_i) g_i$.

\begin{itemize}
\item [2.] For given $\Theta_0$-frame $\seqgr[g]$ for $X_0$
and given sequence $\{X_s\}_{s\in\mn_0}$ (resp.
$\{\Theta_s\}_{s\in\mn_0}$), construct $\{\Theta_s\}_{s\in\mn}$
(resp. $\{X_s\}_{s\in\mn}$) so that $\seqgr[g]$ is a pre-$F$-frame
or $F$-frame for $X_F$ with respect to $\Theta_F$.
\end{itemize}
Let $\seqgr[g]$ be a $\Theta_0$-frame for $X_0$ and let
$\{(X_s,||\cdot||_s)\}_{s\in\mn_0}$ be a sequence of Banach spaces
such that $X_s\subset X_{s-1}\subset X_0$, $s\newin\mn$, in the
topological sense also. The question is whether there exists a
decreasing sequence of Banach sequence spaces
$\{\Theta_s\}_{s\in\mn_0}$ so that $\seqgr[g]$ is a
$\Theta_s$-frame for $X_s$, $s\newin\mn_0$. Also, if
$\Theta_s\subset \Theta_{s-1}\subset \Theta_0$, $s\newin\mn,$ is
given, whether there exists a corresponding sequence
$\{X_s\}_{s\in\mn_0}$ so that $\seqgr[g]$ is a $\Theta_s$-frame
for $X_s$, $s\newin\mn_0$.
 These  problems are interesting in itself
 for Banach spaces,
but if one also imposes the condition that
$\Theta_F=\cap_{s\in\mn_0}\Theta_s$ and $X_F=\cap_{s\in\mn_0} X_s$
are dense in $\Theta_s$, respectively, $X_s$, $s\newin\mn_0,$ then
one comes to pre-$F$- and $F$-frames and to the frame theory for
spaces of test functions and their duals, various spaces of
generalized functions.

The paper is organized as follows. In Section \ref{s2p} we give
the definitions of pre-$F$-  and $F$-  frames for a Fr\'echet
space $X_F$ with respect to a Fr\'echet sequence space $\Theta_F$.
Moreover, we introduce pre-$DF$- and $DF$-frames motivated by the
properties of  dual frames in Banach spaces. The frame expansions
in Fr\'echet spaces and their duals are the subject of Section 3.
We note that for a Fr\'echet space $X_F$, which is not a Banach
space,
 the existence of a pre-$F$-frame $\seqgr[g]$ implies that the corresponding sequence space
$\Theta_F$ must be a Fr\'echet space, which is not a Banach space
(Remark \ref{motiv}).
 Let $\seqgr[g]$ be a $\Theta_0$-frame for $X_0$. In Section \ref{constrfx} we construct
a sequence $\{X_s\}_{s\in\mn}$ for given
$\{\Theta_s\}_{s\in\mn_0}$ and a sequence $\{\Theta_s\}_{s\in\mn}$
for given  $\{X_s\}_{s\in\mn_0}$ so that $\seqgr[g]$ is a
pre-$F$-frame or $F$-frame for $X_F=\cap_{s\in\mn_0}X_s$ with
respect to $\Theta_F=\cap_{s\in\mn_0}\Theta_s$. In Subsection
\ref{subs43} we construct a new sequence space
$\widetilde{\Theta}$, so that given $\Theta$-frame for $X$ (resp.
Banach frame for $X$ with respect to $\Theta$)
 is a $\widetilde{\Theta}$-frame for $X$
(resp. Banach frame for $X$ with respect to $\widetilde{\Theta}$).
With this construction, if $\{X_s\}_{s\in\mn_0}$ is a decreasing
sequence of Banach spaces, we formulate conditions which imply the
existence of a sequence $\{\Theta_s\}_{s\in\mn_0}$ so that
 given $\Theta_0$-Bessel sequence
for $X_0$ is a pre-$F$-frame for $X_F$ with respect to $\Theta_F$
and open questions related to Fr\'echet frames for such sequences
$\{X_s\}_{s\in\mn_0}$ and $\{\Theta_s\}_{s\in\mn_0}$. Our analysis
is illustrated by Propositions \ref{primer1} and \ref{hilbcor},
which show that the assumptions of Theorems \ref{constructone} and
\ref{constructthetas} are essential in the construction of frames.

\section{Preliminaries}\label{s2p}

We use usual notation: $(X, \|\cdot \|)$ is a Banach space and
$(X^*, \|\cdot \|^*)$ is its dual, $( \Theta, \snorm[\cdot]) $ is
a Banach sequence space and  $( \Theta^*, \snorm[\cdot]^*) $ is
the dual of $\Theta$. Recall that $\Theta$ is called {\it solid}
if the conditions $\{c_i\}_{i=1}^\infty\in \Theta$ and $|d_i| \leq
|c_i|,$ $i\newin\mn$, imply that $\{d_i\}_{i=1}^\infty \in \Theta$
and $\snorm[\{d_i\}_{i=1}^\infty] \leq
\snorm[\{c_i\}_{i=1}^\infty].$ If the coordinate functionals on
$\Theta$ are continuous, then $\Theta$ is called a {\it
$BK$-space}. A $BK$-space, for which the canonical vectors form a
Schauder basis, is called a {\it $CB$-space}. If otherwise is not
written, $e_i$ will denote the $i$-th canonical vector
$\{\delta_{ki}\}_{k=1}^\infty$, $i\newin\mn$. A $BK$-space
$\Theta$, which contains all the canonical vectors and for which
there exists $\lambda\geq 1$ such that
$$
\snorm[\{c_i\}_{i=1}^n]_\Theta\leq
\lambda\snorm[\{c_i\}_{i=1}^\infty]_\Theta, \
 n\in\mn, \, \{c_i\}_{i=1}^\infty\in \Theta,
$$ ($\{c_i\}_{i=1}^n\equiv\sum_{i=1}^n c_i e_i$), is called a {\it
$\lambda$--$BK$-space}.

When $\Theta$ is a $CB$-space, then the canonical basis constant
$\sup_{N\in\mn} \|S_N\|$ is a finite number $\geq 1$, where
$$S_N:\Theta\to\Theta,\; S_N(\seqgr[c])=\sum_{i=1}^N c_i e_i, \ N\in\mn $$
(see \cite{Heil}). Every $CB$-space is a $\lambda$--$BK$-space,
where $\lambda$ is the canonical basis constant. If $\Theta$ is a
$CB$-space, then the space $\Theta^\circledast :=\{
\{g(e_i)\}_{i=1}^{\infty} : g\in \Theta^* \}$ with the norm
$\snorm[\{g(e_i)\}_{i=1}^{\infty}]_{\Theta^\circledast}:=\|g\|_{\Theta^*}$
is a $BK$-space, isometrically isomorphic to $\Theta^*$ (see
\cite[p.\,201]{KA}). \label{bkthetastar} When $\Theta$ is a
reflexive $CB$-space, then the coefficient functionals $E_i$,
$i\newin\mn$, associated to the canonical basis $e_i$,
$i\newin\mn$, form a Schauder basis for $\Theta^*$ (see
\cite[p.\,57]{Heil}) and thus $\Theta^\circledast$ is a
$CB$-space, since the canonical vectors
$\{E_j(e_i)\}_{i=1}^\infty$, $j\newin\mn$, form a Schauder basis
for $\Theta^\circledast$. From now on when $\Theta$ is a
$CB$-space, we will always identify $\Theta^*$ with
$\Theta^\circledast$. In the sequel, linear mappings will be
called operators. Recall that an operator $\mathcal{P}:X\to X$ is
called {\it projestion} if $\mathcal{P}^2=\mathcal{P}$.

\subsection*{Pre-Fr\'{e}chet and Fr\'{e}chet frames}

Let $\{Y_s, | \cdot |_s\}_{s\in\mn_0}$ be a sequence of separable
Banach spaces such that \be \label{fx1} \{\nulel\} \neq
\sech[Y]\subseteq \ldots \subseteq Y_2 \subseteq Y_1 \subseteq Y_0
\ee \be  \label{fx2} |\cdot|_0\leq | \cdot |_1\leq | \cdot |_2\leq
\ldots \ee \be \label{fx3} Y_F :=\sech[Y] \;\; \mbox{is dense in}
\;\; Y_s, \;\; s\in\mn_0. \ee

Then $Y_F$ is a Fr\'echet space with the sequence of norms $ |
\cdot |_s, $ $ s\in\mn_0.$
 We will use such sequences in two cases:

1. $Y_s=X_s$ with norm $\|\cdot\|_s, s\in\mn_0;$

2. $Y_s=\Theta_s$ with norm $\snorm[\cdot]_s, s\in\mn_0$.

Note that if the norms $|\cdot|_s$, $s\in\mn_0$, satisfy
 $|\cdot|_s\leq C_s |\cdot|_{s+1}$, $C_s>0$,
 $s\in\mn_0$, then new norms can be
introduced in $Y_s$, $s\newin\mn,$ so that (\ref{fx2}) holds.

\begin{rem} \label{novr}
1. If $\{\Theta_s, \snorm[\cdot]_s\}_{s\in\mn_0}$ is a sequence of
$CB$-spaces, which satisfies (\ref{fx1}) and (\ref{fx2}),
\label{cbspaces} then all finite sequences belong to $\Theta_F$
and form a dense subset of $\Theta_s$; thus (\ref{fx3}) is
automatically satisfied. In this case the canonical vectors form a
Schauder basis of $\Theta_F$, every $\seqgr[c]\in\Theta_F$ can be
written uniquely as $\seqgr[c]=\sumii c_i e_i$
 with the
convergence in $\Theta_s$ for every $s\newin\mn_0$.

2. We refer to \cite[p. 326, pp. 331--332]{MV} for K\" othe type
sequence spaces $\lambda^p(\ma)$ and their duals as a good source
of examples of sequence spaces.
\end{rem}

\begin{defn}
Let $\{X_s, \|\cdot\|_s\}_{s\in\mn_0}$ and $\{\Theta_s,
\snorm[\cdot]_s\}_{s\in\mn_0}$ be sequences of Banach spaces,
which satisfy (\ref{fx1})-(\ref{fx3}). For fixed $s\newin\mn_0$,
an operator $V:\Theta_F\to X_F$ will be called $s$-bounded, if
there exists a constant $K_s>0$ such that $\|V\seqgr[c]\|_s\leq
K_s \snorm[\{c_i\}_{i=1}^\infty]_s$ for all $\seqgr[c]\newin
\Theta_F$. If $V$ is $s$-bounded for every $s\newin\mn_0$, then
$V$ will be called $F$-bounded.
\end{defn}

Note that if $V:\Theta_F\to X_F$ is $F$-bounded, then $V$ is
continuous.
The converse does not hold in general.

 Let $\Theta$ be a $BK$-space. Recall that
$\seqgr[g]\in (X^*)^\mn$ is called a {\it $\Theta$-frame for $X$
with lower bound $A>0$ and upper bound $B<\infty$}, if
\begin{equation} \label{cl}
\{g_i(f)\}_{i=1}^\infty\in\Theta \ \mbox{and} \  A \|f\|\leq
\snorm[\{g_i(f)\}_{i=1}^\infty]\leq B\|f\|, \ f\in X.
\end{equation}

A $\Theta$-frame $\seqgr[g]$ for $X$ is called a {\it Banach frame
for $X$ with respect to $\Theta$} if there exists a continuous
operator $V:\Theta\to X$ such that $V(\{g_i(f)\}_{i=1}^\infty)=f$
for every $f\in X$.

If at least the upper inequality in (\ref{cl}) holds, then
$\seqgr[g]$ is called a {\it $\Theta$-Bessel sequence for $X$ with
 bound $B$}.

We generalize the above definitions to Fr\'echet spaces as
follows:
\begin{defn}\label{fframe}
Let $\{X_s, \|\cdot\|_s\}_{s\in\mn_0}$ be a sequence of Banach
spaces, which satisfies (\ref{fx1})-(\ref{fx3}), and let
$\{\Theta_s, \snorm[\cdot]_s\}_{s\in\mn_0}$ be a sequence of
$BK$-spaces, which satisfies (\ref{fx1})-(\ref{fx3}). A sequence
$\seqgr[g]\newin ({X_F^*})^\mn$ is called a pre-$F$-frame for
$X_F$ with respect to $\Theta_F$ if for every $s\in\mn_0$ there
exist constants $0<A_s\leq B_s<\infty$ such that
\begin{equation}\label{fframestar}
\{g_i(f)\}_{i=1}^\infty\in\Theta_F, \ f\newin X_F,
\end{equation}
\begin{equation}\label{fframetwostar}
A_s \|f\|_s\leq \snorm[\{g_i(f)\}_{i=1}^\infty]_{s}\leq
B_s\|f\|_s, \ f\newin X_F.
\end{equation}
The constants $B_s$ (resp. $A_s$), $s\in\mn_0$, are called upper
(resp. lower) bounds for $\seqgr[g]$. A pre-$F$-frame is called
{\it tight}, if $A_s=B_s, s\newin\mn_0$.

Moreover, if there exists an $F$-bounded operator
$V:\Theta_F\rightarrow X_F$ so that $V(\{g_i(f)\}_{i=1}^\infty)=f$
for all $f\newin X_F,$ then a pre-$F$-frame $\seqgr[g]$ is called
an $F$-frame (Fr\'echet frame) for $X_F$ with respect to
$\Theta_F$ and $V$ is called an $F$-frame operator for
$\seqgr[g]$.

When (\ref{fframestar}) and at least the upper inequality in
(\ref{fframetwostar}) hold, then $\seqgr[g]$ is called an
$F$-Bessel sequence for $X_F$ with respect to $\Theta_F$ with
bounds $B_s$, $s\newin\mn_0$.
\end{defn}

If $ X = X_F = X_s $ and $ \Theta =  \Theta_F =  \Theta_s, $ $ s
\in \mn_0,$ the above definitions of a pre-$F$-frame, $F$-frame
and $F$-Bessel sequence for $X_F$ give the definitions of a
$\Theta$-frame, Banach frame and $\Theta$-Bessel sequence for $X$,
respectively.

Note that if $\Theta$ is a reflexive Banach space and $\seqgr[g]$
is a $\Theta$-frame for a Banach space $X$, then $X$ is also
reflexive because it is isomorphic to a closed subspace of
$\Theta$. \label{refls}

Let $\{X_s\}_{s\in\mn_0}$ be a sequence of Banach spaces which
satisfies (\ref{fx1})-(\ref{fx3}), $\{\Theta_s\}_{s\in\mn_0}$ be a
sequence of $BK$-spaces which satisfies (\ref{fx1})-(\ref{fx3})
and $\seqgr[g]\newin(X_F^*)^\mn$ be an $F$-Bessel sequence for
$X_F$ with respect to $\Theta_F$. Fix $i\newin\mn$ and
$s\newin\mn_0$. Since $\Theta_s$ is a $BK$-space, the $i$-th
coordinate functional on $\Theta_s$ is bounded and thus there
exists $K_{i,s}>0$ such that $|g_i(f)|\leq
K_{i,s}\snorm[\{g_i(f)\}_{i=1}^\infty]_s$, $f\newin X_F$.
Therefore, by the upper inequality in (\ref{fframetwostar}), $g_i$
is bounded on $X_F$ with respect to the norm $\|\cdot\|_s$.
\label{conts} Now by (\ref{fx3}) it follows that $g_i$ has a
unique continuous extension on $X_s$ which will be denoted by
$g^s_i.$ Thus, for every $s\newin\mn_0$, $g_i^s\newin X_s^*$,
$i\newin\mn$, and moreover, $g_i^{s}|_{X_{{t}}}=g_i^{t}$,
$i\newin\mn$, for $t>s$.

\begin{prop} {\rm \cite{pst}} \label{gs}
Let $\{X_s\}_{s\in\mn_0}$ be a sequence of Banach spaces, which
satisfies (\ref{fx1})-(\ref{fx3}), and let
$\{\Theta_s\}_{s\in\mn_0}$ be a sequence of
$\lambda_s$--$BK$-spaces, which satisfies (\ref{fx1})-(\ref{fx3}).
If $\seqgr[g]\newin (X_F^*)^\mn$ is an $F$-Bessel sequence (resp.
pre-$F$-frame) for $X_F$ with respect to $\Theta_F$ with bounds
$B_s$ (resp. lower bounds $A_s$ and upper bounds $B_s$),
$s\newin\mn_0$, then for any $s\newin\mn_0$,
$\seq[g^s]_{i=1}^\infty$ is a $\Theta_s$-Bessel sequence (resp.
$\Theta_s$-frame) for $X_s$ with bound $\lambda_s B_s$ (resp.
lower bound $A_s$ and upper bound $\lambda_s B_s$).
 \end{prop}

In the case of a single Banach space ($ X = X_F = X_s $, $ \Theta
=  \Theta_F =  \Theta_s$, $s\newin\mn_0$) the following assertion
holds.

\begin{prop}
{\rm\cite{Sduals}}\label{bnew} Let $\Theta$ be a $CB$-space and
$\seqgr[g]\newin {(X^*)}^\mn$ be a $\Theta$-Bessel sequence for
$X$. Consider the  conditions:
\begin{itemize}
\item[$(\mathcal{P}_1)$]  $\seqgr[g]$ is a Banach frame for X with
respect to $\Theta$.
\item[$(\mathcal{P}_2)$] There exists a $\Theta^*$-Bessel sequence $\seqgr[f]$  ($f_i\newin X\subseteq X^{**}, i\in\mn$) for $X^*$ such
that
$
f=\sumii g_i(f)f_i, \ \  f \in X.
$
\item[$(\mathcal{P}_3)$] There exists a $\Theta^*$-Bessel sequence $\seqgr[f]$ ($f_i\newin X
\subseteq X^{**}, i\in\mn$) for $X^*$
 such that
$
g=\sumii g(f_i)g_i, \ \ g \in X^*.
$
\end{itemize}
Then:
\begin{itemize}
\item[$(a)$]  $(\mathcal{P}_1)\Leftrightarrow (\mathcal{P}_2)$.
\item[$(b)$] If both $\Theta$ and $\Theta^*$ are $CB$-spaces,
then $(\mathcal{P}_1)\Leftrightarrow (\mathcal{P}_2)
\Leftrightarrow (\mathcal{P}_3)$ and each one of the conditions
$(\mathcal{P}_2)$ and $(\mathcal{P}_3)$ implies that $\seqgr[f]$
is a $\Theta^*$-frame for $X^*$, called {\it a dual of the
$\Theta$-frame $\seqgr[g]$}.
\item[$(c)$]
If $\Theta$ is a reflexive $CB$-space,  then each one of the
conditions $(\mathcal{P}_2)$ and $(\mathcal{P}_3)$ implies that
$\seqgr[f]$ is a Banach frame for $X^*$ with respect to
$\Theta^*$.
\end{itemize}
\end{prop}

We will extend the above proposition to Fr\'{e}chet spaces and
obtain series expansions by the use of a Fr\'{e}chet frame and the
corresponding  dual frame.

Once again we note that the sequence $\seqgr[f]$ in Proposition
\ref{bnew} has dual properties in comparison to $\seqgr[g]$: $f_i$
belongs to $X$ while $g_i$ belongs to $X^*$, $i\newin\mn$,
$\seqgr[f]$ is a $\Theta^*$-Bessel sequence for $X^*$ while
$\seqgr[g]$ is a $\Theta$-Bessel sequence for $X$. Having this in
mind as well as Proposition \ref{gs}, we give the following
definition:

\begin{defn}\label{fframedual}
Let all the assumptions of Definition \ref{fframe} hold and assume
moreover, that $\Theta_s$, $s\newin\mn_0$, are $CB$-spaces (then
their duals are $BK$-spaces). A sequence $\seqgr[f]\newin
({X_F})^\mn$ is called a

$DF$-Bessel sequence for $X_F^*$ with respect to $\Theta_F^*$ if
it is a $\Theta_s^*$-Bessel sequence for $X_s^*$ for every
$s\in\mn_0$;

 pre-$DF$-frame for $X_F^*$ with respect to $\Theta_F^*$ if
it is a $\Theta_s^*$-frame for $X_s^*$ for every $s\in\mn_0$;

 $DF$-frame for $X_F^*$ with respect to $\Theta_F^*$ if it is a Banach frame for $X_s^*$ with respect to $\Theta_s^*$ for every $s\in\mn_0$.
\end{defn}

\section{Series Expansions} \label{se}

Throughout this section we assume:
\begin{itemize}
\item [1.]  $\{X_s,
\|\cdot\|_s\}_{s\in\mn_0}$ is a sequence of Banach spaces,
satisfying (\ref{fx1})-(\ref{fx3});
\item[2.] $\{\Theta_s, \snorm[\cdot]_s\}_{s\in\mn_0}$
is a sequence of $\lambda_s$--$BK$-spaces, satisfying
(\ref{fx1})-(\ref{fx3});
\item[3.] $\seqgr[g]\newin (X_F^*)^\mn$.
\end{itemize}

Let $\seqgr[g]$ be a pre-$F$-frame for $X_F$ with respect to
$\Theta_F$ with lower bounds $A_s$ and upper bounds $B_s$,
$s\newin\mn_0$. By Proposition \ref{gs}, for every $s\in\mn_0$,
\begin{equation} \label{asb2}
A_s \|f\|_s\leq \snorm[\{g_i^s(f)\}_{i=1}^\infty]_{s}\leq
\lambda_s B_s\|f\|_s, \ \
 f\in X_{s}.
\end{equation}
By (\ref{asb2}) and (\ref{fframetwostar}), one can consider the
following continuous operators
\begin{eqnarray}
U_s : X_s \to \Theta_s,& & U_sf=\{g_i^s(f)\}_{i=1}^\infty, \ s\in\mn_0,\label{operatorus}\\
U: X_F \to \Theta_F, & & Uf=\{g_i(f)\}_{i=1}^\infty,
\label{operatoru}
\end{eqnarray}
and their inverses
$$
U_s^{-1} : \rus \to X_s, \ \ U^{-1}: \ru\to X_F.
$$
For every $s\newin\mn_0$, the operator $U_s^{-1}$ is bounded with
$\|U_s^{-1}\|\leq \frac{1}{A_s}$.

\begin{lem} \label{ruclosed}
Let $\seqgr[g]$ be a pre-$F$-frame for $X_F$ with respect to
$\Theta_F$.
 The range $\ru$ of the operator \,$U$, defined by
(\ref{operatoru}), is closed in $\Theta_F$ and the inverse
operator $U^{-1}$ is $F$-bounded.
\end{lem}
\bp Denote the lower bounds for $\seqgr[g]$ by $A_s,
s\newin\mn_0$, and the upper bounds by $B_s, s\newin\mn_0$. Let
$f_n\newin X_F$, $n\newin\mn$, be such that
$\{Uf_n\}_{n=1}^\infty$  converges in $\Theta_F$ as $n\to\infty$.
Fix $\varepsilon>0$ and $s\in\mn_0$. Then there exists
$N_0(s)\newin\mn$ such that
$\snorm[U_sf_n-U_sf_k]_{s}<\varepsilon$ for every $n,k > N_0 (s)$.
By the lower inequality in (\ref{asb2}), we obtain $$0\leq A_s
\|f_n-f_k\|_s <\varepsilon, \ \ n,k >N_0(s)$$ and hence
$\{f_n\}_{n=1}^\infty$ converges in $X_s$ as $n\to\infty$.
Therefore, $\{f_n\}_{n=1}^\infty$ converges to some $f$ in $X_F$.
Hence, $Uf\in\ru$ is the limit of $\{Uf_n\}_{n=1}^\infty$ in
$\Theta_F$. This implies that $\ru$ is closed in $\Theta_F$.

Since for every  $\{c_i\}_{i=1}^\infty\in \ru$ and every
$s\newin\mn_0$,
$$ \|U^{-1}\seqgr[c]\|_s = \|U_s^{-1}\seqgr[c]\|_s\leq 1/A_s \snorm[\{c_i\}_{i=1}^\infty]_{s},$$
it follows that $U^{-1}$ is $F$-bounded.
 \ep

\begin{rem} \label{motiv}
Assume that $\Theta_s=\Theta$ for every $s\in\mn_0$ and that
$\seqgr[g]$ is a pre-$F$-frame for $X_F$ with respect to
$\Theta_F=\Theta$. In this case $X_F$ is isomorphic to $\ru$,
which is a closed subspace of $\Theta$ by the above lemma. Hence,
$X_F$ is isomorphic to a Banach space. This explains that in the
case when $X_F$ is not a Banach space, $\Theta_F$ must not be a
Banach space. Thus, different sequence spaces $\Theta_s$ must be
used for the construction of $\Theta_F$.
\end{rem}

\begin{rem}
\label{ffrp} Consider a pre-$F$-frame $\seqgr[g]$ for $X_F$ with
respect to
 $\Theta_F$. The mapping $V$ (if it exists) in Definition \ref{fframe} is
an $F$-bounded extension of\, $U^{-1}$ and the operator $UV$ is
 an $F$-bounded projection from $\Theta_F$ onto $\ru$.
Conversely,
 assume that there exists an $F$-bounded projection $\mathcal{P}$ from $\Theta_F$ onto $\ru$. In this case
 the operator $U^{-1} \mathcal{P}:\Theta_F\to X_F$ is an $F$-bounded extension of $U^{-1}:\ru\to X_F$.
 Therefore, a pre-$F$-frame for $X_F$ with respect to $\Theta_F$ is an $F$-frame for $X_F$ with respect to $\Theta_F$ if and only if
there exists an $F$-bounded projection from $\Theta_F$ onto $\ru$.
\end{rem}

\begin{thm} \label{diff} Let
$\seqgr[g]$ be an $F$-frame for $X_F$ with respect to $\Theta_F$.
Then the following holds.
\begin{itemize}
\item[{\rm(a)}] For every $s\newin\mn_0$, the sequence $\{g_i^s\}_{i=1}^\infty$
is a Banach frame for $X_s$ with respect to $\Theta_s$.
\item[{\rm(b)}] If \,$\Theta_s$, $s\newin\mn_0$,
are  $CB$-spaces, then there exists  $\seqgr[f]\newin (X_F)^\mn$,
which is a $DF$-Bessel sequence for $X_F^*$ with respect to
$\Theta_F^*$
 such that
\begin{eqnarray}
f&=&\sum_{i=1}^{\infty} g_i(f) f_i,  \   f\newin X_F, \
\mbox{(in $X_F$),} \label{frepr} \\
g&=&\sum_{i=1}^{\infty} g(f_i) g_i,  \  g\in X_F^*, \
\mbox{(in $X_F^*$),} \label{frepr2}\\
f&=&\sum_{i=1}^{\infty} g_i^s(f) f_i, \  f\newin X_s, \
 s\newin\mn_0.\label{fsrepr}
\end{eqnarray}
\item[{\rm(c)}] If \,$\Theta_s$ and  $\Theta^*_s$, $s\newin\mn_0$, are  $CB$-spaces,
then there exists $\seqgr[f]\newin (X_F)^\mn$, which is a
pre-$DF$-frame for $X_F^*$ with respect to $\Theta_F^*$ such that
(\ref{frepr})-(\ref{fsrepr}) hold and moreover,
\begin{equation} \label{xs}
g=\sum_{i=1}^{\infty} g(f_i) g_i^s, \ g\in X_s^*,\,
 s\in\mn_0. \ \ \ \ \ \ \ \ \ \ \ \ \ \ \ \ \ \ \ \ \ \ \ \ \ \ \
\end{equation}
\item[{\rm(d)}]
If \,$\Theta_s$, $s\newin\mn_0$, are reflexive $CB$-spaces, then
there exists $\seqgr[f]\newin (X_F)^\mn$, which is a $DF$-frame
for $X_F^*$ with respect to $\Theta_F^*$ such that
(\ref{frepr})-(\ref{xs}) hold.
\end{itemize}
\end{thm}

\bp Having in mind Remark \ref{ffrp}, let $\mathcal{P}$ denote an
$F$-bounded projection from $\Theta_F$ onto $\ru$. Fix
$s\in\mn_0$. By Proposition \ref{gs},  $\{g_i^s\}_{i=1}^\infty$ is
a $\Theta_s$-frame for $X_s$ and thus we can consider the operator
$U_s$, given by (\ref{operatorus}). Consider $\mathcal{P}$ as an
operator from $\Theta_F\subseteq \Theta_s$ into $R(U_s)$ and note
that $R(U_s)$ is a closed subspace of $\Theta_s$. Since $\Theta_F$
is dense in $\Theta_s$ and $\ru$ is dense in $R(U_s)$,
$\mathcal{P}$ has unique continuous extension $\mathcal{P}_s$
defined from $\Theta_s$ into $R(U_s)$ and moreover,
$\mathcal{P}_s$ is a projection from $\Theta_s$ onto $R(U_s)$.
Define $V_s:\Theta_s\to X_s$ by $V_s:=U_s^{-1}\mathcal{P}_s$.

(a) For every $s\newin\mn_0$, $V_s$ is a continuous extension of
$U_s^{-1}$ and hence the $\Theta_s$-frame $\{g_i^s\}_{i=1}^\infty$
is a Banach frame for $X_s$ with respect to $\Theta_s$.

(b) Consider the  operator $V=\uinv \mathcal{P}$, which is an
$F$-bounded extension of $U^{-1}$. Since all the canonical vectors
belong to $\Theta_F$, for $i\newin \mn$ we define $f_i:=V e_i$.
Let $f\in X_F$. By Remark \ref{novr} applied to $\seqgr[c]=Uf
\newin\Theta_F$,
$$\sum_{i=1}^n
g_i(f) e_i\to Uf \ \, \mbox{in} \, \ \Theta_F \ \, \mbox{as} \, \
n\to\infty.$$ The continuity of $V$ implies that
$$V\left(\sum_{i=1}^n g_i(f)e_i\right)\to VUf=f \ \mbox{in}\ X_F\ \mbox{as}\ n\to\infty.$$
Hence, $\sum_{i=1}^n g_i(f) f_i\to f$ in $X_F$, as $n\rightarrow
\infty$. Moreover, (\ref{frepr}) implies (\ref{frepr2}).

Let $s\newin\mn_0$.
 Since $\mathcal{P}_s$ is an extension of $\mathcal{P}$ and  $U_s^{-1}$ is an extension of $\uinv$, it follows that
$V_s$ is an extension of $V$. Therefore, $V_s e_i=Ve_i=f_i$, $
i\newin\mn$.
 By the fact that $\seqgr[e]$ is a
Schauder basis for $\Theta_s$, it follows that for every $f\newin
X_s$,
$$f=V_sU_s f=
V_s\left(\sumii g_i^s(f)e_i\right)=\sum_{i=1}^{\infty}
g_i^s(f)V_se_i= \sum_{i=1}^{\infty} g_i^s(f)f_i.$$ Let now $g\in
X_s^*$. Then $gV_s\in\Theta_s^*$, $\{g(f_i)\}_{i=1}^\infty=
\{(gV_s) e_i\}_{i=1}^\infty\newin\Theta_s^*$ and
$$\snorm[\{g(f_i)\}_{i=1}^\infty]^*_s= \snorm[gV_s]^*_s \leq
\|g\|_s^*\, \,\|V_s\|.$$ Therefore, $\seqgr[f]$ is a
$\Theta_s^*$-Bessel sequence for $X_s^*$.

(c) Let $s\newin\mn_0$.  Proposition \ref{bnew}(b) implies that
the sequence $\seqgr[f]$, determined in (b), is a
$\Theta_s^*$-frame for $X_s^*$. Now (\ref{xs}) follows from
(\ref{fsrepr}) and \cite[Lemma 4.3]{Sduals}.

(d) Consider the sequence $\seqgr[f]$ defined in (b). Let $s\newin
\mn_0$. By Proposition \ref{bnew}(c) it follows that $\seqgr[f]$
is a Banach frame for $X_s^*$ with respect to $\Theta_s^*$. \ep

Note, for validity of (\ref{frepr}) and (\ref{frepr2}) it is
enough to
 assume that there exits a continuous projection $\mathcal{P}$  from $\Theta_F$
 onto $\ru$. The $F$-boundedness of $\mathcal{P}$ is essential for
 (\ref{fsrepr}), (\ref{xs}) and (a).
\begin{rem}
If $\seqgr[g]$ is an $F$-frame for $X_F$ with respect to
$\Theta_F$ and $\Theta_s$, $s\newin\mn_0$, are $CB$-spaces,
Theorem \ref{diff} (a) and Proposition \ref{bnew} imply that for
every $s\newin\mn_0$ there exists a sequence $\seqgr[f^s]\newin
X_s^\mn$ such that $f=\sumii g_i^s(f)f_i^s$ for every $f\newin
X_s$. The sequences $\seqgr[f^s]$ might be different  for
different $s$. Since our aim is to obtain series expansions in
$X_F$, we need a same sequence $\seqgr[f]$ satisfying $f=\sumii
g_i(f) f_i$, $f\newin X_F$, with convergence in $\|\cdot\|_s$-norm
for every $s\newin\mn_0$. Because of that we have Theorem
\ref{diff} (b) determining a sequence $\seqgr[f]\newin X_F^\mn$,
which gives series expansions on every level $s$, $s\newin\mn_0$.
\end{rem}

In order to have series expansions in $X_F$ via a pre-$F$-frame
(or $F$-Bessel sequence) $\seqgr[g]$ and a $DF$-Bessel sequence
$\seqgr[f]$, one must have that $\seqgr[g]$ is an $F$-frame. This
will be proved in the next theorem.

\begin{thm}\label{nec} Let $\Theta_s$, $s\newin\mn_0$, be $CB$-spaces
and  $\seqgr[g]$ be an $F$-Bessel sequence for $X_F$ with respect
to $\Theta_F$. Then the following holds.
\begin{itemize}
\item[{\rm (a)}] There exists $\seqgr[f]\newin (X_F)^\mn$, such that it
is a $DF$-Bessel sequence for $X_F^*$ with respect to $\Theta_F^*$
and satisfies (\ref{frepr})
 if and only if
 $\seqgr[g]$ is an $F$-frame for $X_F$ with respect to
$\Theta_F$.
\item[{\rm (b)}] Let $\seqgr[f]\newin (X_F)^\mn$ be a $DF$-Bessel
sequence for $X_F^*$ with respect to $\Theta_F^*$ which satisfies
(\ref{frepr}). If both $\Theta_s$ and $\Theta_s^*$ are
$CB$-spaces, $s\in\mn_0$, (resp. $\Theta_s$ is a reflexive
$CB$-space, $s\in\mn_0$), then $\seqgr[f]$ is a pre-$DF$-frame
(resp. $DF$-frame) for $X_F^*$ with respect to $\Theta_F^*$.
\end{itemize}
\end{thm}
\bp (a) If  $\seqgr[g]$ is an $F$-frame for $X_F$ with respect to
$\Theta_F$, the assertion is given by Theorem \ref{diff}.

For the converse, assume that $\seqgr[f]\newin (X_F)^\mn$ is a
$DF$-Bessel sequence for $X_F^*$ with respect to $\Theta_F^*$
which satisfies (\ref{frepr}).
 First we will prove that (\ref{fsrepr}) holds. Fix $s\newin \mn_0$. By
Proposition \ref{gs}, $\{g^s_i\}_{i=1}^{\infty}$ is a
$\Theta_s$-Bessel sequence for $X_s$ and thus the operator $U_s$,
given by (\ref{operatorus}), is bounded. Let $\varepsilon>0$ and
$f\newin X_s$. Since $X_F$ is dense in $X_s$, there exists
$h\newin X_F$ such that $\|f-h\|_s<\varepsilon$. Applying
(\ref{frepr}) to $h$, there exists $N\in\mn$ such that
$$\left\|h-\sum_{i=1}^n g_i(h) f_i\right\|_s< \varepsilon, \  n>N.$$
Let $\lambda_s$ denote the canonical basis constant in $\Theta_s$
(see the preliminaries). Since
$\{g_i^s(h-f)\}_{i=1}^\infty\newin\Theta_s$, we have
\begin{eqnarray*}
\snormb[\sum_{i=1}^n g_i^s(h-f) e_i]_{s} &\leq& \lambda_s
\snormb[\sum_{i=1}^\infty g_i^s(h-f) e_i]_{s}=\lambda_s
\snormb[U_s(h-f)]_{s}\\
&<&  \varepsilon \, \lambda_s\|U_s\|.
\end{eqnarray*}
Since $\seqgr[f]$ is a $\Theta_s^*$-Bessel sequence for $X_s^*$,
by \cite[Proposition 3.2]{CCS}, 
the operator $T_s:\seqgr[d]\to\sumii d_i f_i$ is bounded from
$\Theta_s$ into $X_s$. Thus for $n>N$,
\begin{eqnarray*}
\left\|f-\sum_{i=1}^n g_i^s(f) f_i\right\|_s &\leq& \|f-h\|_s +
\left\|h-\sum_{i=1}^n g_i(h) f_i\right\|_s +
 \left\|\sum_{i=1}^n
g_i^s(h-f) f_i\right\|_s\\
&\leq& 2\varepsilon + \left\|T_s\left(\sum_{i=1}^n g_i^s(h-f)
e_i\right)
\right\|_s \\
&<& \varepsilon (2+ \lambda_s\|T_s\|\,\|U_s\|)
\end{eqnarray*}
and hence $f=\sumii g_i^s(f) f_i$. Now Proposition \ref{bnew}(a)
implies that $\{g_i^s\}_{i=1}^\infty$ is a Banach frame for $X_s$
with respect to $\Theta_s$.

In order to conclude that $\seqgr[g]$ is an $F$-frame, it remains
to prove that there exists an $F$-bounded extension of $U^{-1}$.
For every $\seqgr[d]\in\Theta_F$, the series $\sumii d_i f_i$
converges in $X_s$ for every $s\in\mn_0$ and therefore it
converges in $X_F$. Thus the operator $T$, given by $T
(\seqgr[d]):=\sumii d_i f_i$, is defined  from $\Theta_F$ into
$X_F$. For every $s\in\mn_0$, $T=T_s\vert_{\Theta_F}$ and hence
$T$ is $s$-bounded.
 Since $\seqgr[e]$ is a Schauder basis for $\Theta_F$ and $T$ is continuous on $\Theta_F$, it
follows that $Uf=\sumii g_i(f) e_i$ and
$$TU f=T\left(\sumii g_i(f) e_i\right)=\sumii g_i(f) Te_i=\sumii g_i(f) f_i=f, \
\,f\in X_F,$$ i.e. $T$ is an extension of $U^{-1}$.

(b) Assume now that both $\Theta_s$ and $\Theta_s^*$ are
$CB$-spaces. By the proof in (a), it follows that (\ref{fsrepr})
holds. Now Proposition \ref{bnew}(b) implies that $\seqgr[f]$ is a
$\Theta_s^*$-frame for $X_s^*$. Moreover, if $\Theta_s$ is a
reflexive $CB$-space, then Proposition \ref{bnew}(c) implies that
$\seqgr[f]$ is a Banach frame for $X_s^*$ with respect to
$\Theta_s^*$. \ep

\v Note that Theorems \ref{diff} and \ref{nec} clarify Theorem 5.3
in \cite{pst}, where this theorem was quoted without proof.

\begin{rem} Throughout this section we have considered $\seqgr[g]\newin
(X_F^*)^\mn $ which is an $F$-Bessel sequence for $X_F$ with
respect to $\Theta_F$, at least. For such a sequence $\seqgr[g]$
we have that $\{g_i^0\}_{i=1}^\infty\newin (X_0^*)^\mn$ and
$g_i^s=g_i^0\vert_{_{X_s}}$. If instead of $\seqgr[g]\newin
(X_F^*)^\mn $ it is assumed that $\seqgr[g]$ belongs to
$(X_0^*)^\mn$, then (\ref{fx2}) implies that $g_i\vert_{X_s}$ is
continuous on $X_s$ for every $s\newin\mn$ and thus
$\{g_i\vert_{X_F}\}_{i=1}^\infty$ belongs to $(X_F^*)^\mn$.
Therefore, all the assertions in the present section can be stated
with $\seqgr[g]\newin (X_0^*)^\mn $ instead of $\seqgr[g]\newin
(X_F^*)^\mn $ and $\{g_i\vert_{_{X_s}}\}_{i=1}^\infty$ instead of
$\seqgr[g^s]$.
\end{rem}

We end this section with an example in which we use notation,
notions and results of \cite{AST}.
\begin{exmp}\label{ast}
{\rm
 Let $\mathcal{L}^\infty:=\left\{f\ :\
\sup_{x\in[0,1]}\sum_{k\in\mz}
 \mid f(x-k) \mid < \infty\right\}.$
When $\phi\in\mathcal{L}^\infty$ and $p\in(1,\infty)$, the series
$\sum_{k\in\mz} c_k \phi(\cdot - k)$ converges for every
$\{c_k\}_{k\in\mz}\newin\ell^p(\mz)$ (in $L^p$-sense) and
\begin{equation}\label{vpspace}
V_p(\phi):=\left\{ \sum_{k\in\mz} c_k \phi(\cdot - k) \ : \
\{c_k\}_{k\in\mz}\in\ell^p(\mz) \right\}
\end{equation}
is a  shift invariant subspace of $L^p(\mr)$, but not necessarily
closed in $L^p(\mr)$. Assertions in \cite{AST} state that} if
$\phi\in\mathcal{L}^\infty$ and $\{\phi(\cdot - k)\}_{k\in\mz}$ is
a $p$-frame for $V_p(\phi)$, then:
\begin{itemize}
\item[1)] $V_p(\phi)$ is closed in $L^p$;
\item[2)] $\{\phi(\cdot - k)\}_{k\in\mz}$ is an
$r$-frame for $V_r(\phi)$ for every $r\in(1,\infty)$;
\item[3)] $\exists$  $\psi\in\mathcal{L}^\infty$
 such that $\{\psi(\cdot - k)\}_{k\in\mz}$ is a $p$-frame for
$V_p(\psi)=V_p(\phi)$ and
$$ f = \sum_{k\in\mz} \<f, \psi(\cdot - k)\> \,\phi(\cdot - k)
= \sum_{k\in\mz} \<f, \phi(\cdot - k)\> \,\psi(\cdot - k), \
\forall f \in V_p(\phi).
$$
\end{itemize}
Moreover, the function $\psi$ belongs to $V_1(\phi)\subset
V_p(\phi)$ and it is independent of $p$, $p\newin (1,\infty)$.

{\rm
 Let now $\phi\in\mathcal{L}^\infty$ and
$\{g_k\}_{k\in\mz}=\{\phi(\cdot - k)\}_{k\in\mz}$ be a $2$-frame
for $V_2(\phi)$. As an example of a sequence of subspaces of
$X_0=V_2(\phi)$ and sequence of subspaces of $\Theta=\ell^2$,
consider
$$X_s=V_{1+1/(s+1)}(\phi) \mbox{ and } \Theta_s=\ell^{1+1/(s+1)}(\mz),\; s\newin\mn_0.$$
By the well known properties of the $\ell^p$-spaces, it follows
that $\Theta_{s+1}\subset \Theta_s$ and $\snorm[\cdot]_s\leq
\snorm[\cdot]_{s+1}$. Moreover, $\Theta_F\supset \ell^1$ and
finite sequences are dense in $\Theta_F$. By 2),
$\{g_k\}_{k\in\mz}$ is a $\Theta_s$-frame for $X_s$ for every
$s\in\mn_0$; let $A_s$ and $B_s$, $s\newin\mn_0$, denote lower and
the upper bounds, respectively. Let $s\newin\mn_0$. By
(\ref{vpspace}) and the inclusion $\ell^{1+1/(s+2)}\subset
\ell^{1+1/(s+1)}$, one obtains that $X_{s+1}\subset X_s$ and
moreover, for every $f\newin X_s$ one has
$$\|f\|_s\leq \frac{1}{A_s} \snorm[\{g_k(f)\}_{k\in\mz}]_{s}\leq \snorm[\{g_k(f)\}_{k\in\mz}]_{s+1}\leq \frac{B_{s+1}}{A_s}\|f\|_{s+1},$$
which implies that $\|\cdot\|_s\leq C_s\|\cdot\|_{s+1}$ for some
constant $C_s>0$, $s\in\mn_0$. Also, the set $\{f=\sum_{j\in F}
c_j \phi(\cdot - j) \,:\, F \ \mbox{is a finite subset of}\ \mz
\}$ is a subset of $X_F$ and it is dense in any $X_s$,
$s\in\mn_0$. Thus, $\{g_k\}_{k\in\mz}$ is a pre-$F$-frame for
$X_F$ with respect to $\Theta_F$, where
$$\Theta_F=\cap_{s\in\mn_0}\ell^{1+1/(s+1)}(\mz) \mbox{ and }
X_F=\cap_{s\in\mn_0}V_{1+1/(s+1)}(\phi).$$ Let
$\{f_k\}_{k\in\mz}:=\{\psi(\cdot - k)\}_{k\in\mz}$, where $\psi$
is given by 3). Then (\ref{frepr}) holds  and Theorem \ref{nec}
implies that $\seqgr[g]$ is an $F$-frame for $X_F$ with respect to
$\Theta_F$.
}
\end{exmp}

\section{Constructions of sequences of spaces}
\label{constrfx}

For the analysis of series expansions, it is reasonable to have
that the sequence space is a $CB$-space, which is not always the
case. But if $\seqgr[g]\newin (X^*)^\mn$ and there exists a
sequence $\seqgr[f]\newin X^\mn$ such that $f=\sum_{i=1}^{\infty}
g_i(f)f_i$, for all $f\newin X$, then there exists a $CB$-space
$\Theta$ such that $\seqgr[g]$ is a Banach frame for $X$ with
respect to $\Theta$ (see \cite[Proposition 2.9]{CCS}). Below we
apply this assertion to the Fr\'echet case.

 \begin{prop}\label{constr1}
Let $\{X_s\}_{s\in\mn_0}$ be a sequence of Banach spaces, which
satisfies (\ref{fx1})-(\ref{fx3}) and $\seqgr[g]\in (X_0^*)^\mn$.
Assume that there is a sequence $\seqgr[f]\in (X_F\setminus
\nulel)^\mn$ such that for every $s\in\mn_0$ and every $f\in X_s$,
\begin{equation}\label{fxf}
f=\sum_{i=1}^\infty g_i(f) f_i \ \,\mbox{in} \,\ X_s.
\end{equation}
Then there exists a sequence $\{\Theta_s\}_{s\in\mn_0}$ of
$CB$-spaces which satisfies (\ref{fx1})-(\ref{fx3}) and such that
$\{g_i\vert_{_{X_F}}\}_{i=1}^\infty$ is an $F$-frame for $X_F$
with respect to $\Theta_F$. Moreover, $\seqgr[f]$ is a $DF$-Bessel
sequence for $X_F^*$ with respect to $\Theta_F^*$.
\end{prop}
\bp For every  $s\newin\mn_0$, put
\begin{equation}\label{ts}
\Theta_s:=\left\{ \seqgr[c] : \left\{\sum_{i=1}^n c_i
f_i\right\}_{n=1}^\infty \ \mbox{converges in}\ X_s\right\},
\end{equation}
$$
\snorm[\{c_i\}_{i=1}^\infty]_s:= \sup_n \left\|\sum_{i=1}^n c_i
f_i\right\|_s.
$$
By \cite[Lemma 3.5 and Proposition 2.9]{CCS} applied to $X_s$, it
follows that $\Theta_s$ is a $CB$-space and
$\{g_i\vert_{_{X_s}}\}_{i=1}^\infty$ is a Banach frame for $X_s$
with respect to $\Theta_s$. Since $X_{s+1}\subseteq X_{s}$ and
$\|\cdot\|_{s}\leq \|\cdot\|_{s+1}$, it follows that
$\Theta_{s+1}\subseteq \Theta_s$ and $\snorm[\cdot]_{s}\leq
\snorm[\cdot]_{s+1}$. By Remark \ref{novr},
$\{\Theta_s\}_{s\in\mn_0}$ satisfies (\ref{fx3}). By (\ref{ts}),
$\sum_{i=1}^n c_i f_i$ converges in $X_s$ as $n\to\infty$ for
every $\seqgr[c]\newin \Theta_s$. Thus, \cite[Proposition
3.2]{CCS}
 implies that $\seqgr[f]$ is a  $\Theta_s^*$--Bessel sequence for $X_s^*$.
 By Theorem \ref{nec}(a), $\seqgr[g]$ is an $F$-frame for $X_F$ with
respect to $\Theta_F$. \ep

\subsection{Construction of $\{X_s\}_{s\in\mn_0}$}
We start with a sequence $\{\Theta_s\}_{s\in\mn_0}$ and a
$\Theta_0$-frame $\seqgr[g]$ for $X_0$ in order to construct
$\{X_s\}_{s\in\mn}$ such that $\seqgr[g]$ is an $F$-frame for
$X_F$ with respect to $\Theta_F$.

\begin{thm} \label{constructx}
Let $X_0\neq\{\nulel\}$ be a Banach space,
$\{\Theta_s\}_{s\in\mn_0}$ be a sequence of $BK$-spaces which
satisfies (\ref{fx1})-(\ref{fx3}) and $\seqgr[g]\newin
(X^*_0)^\mn$ be a $\Theta_0$-frame for $X_0$ with respect to
$\Theta_0$ with bounds $1\leq A_0\leq B_0<\infty$. Assume that
$$M:= \{\{g_i(f)\}_{i=1}^\infty:f\in X_0\} \cap \Theta_F$$ is
dense in $\{\{g_i(f)\}_{i=1}^\infty:f\in X_0\}\cap \Theta_s\neq
\{\nulel\}$ with respect to the $\snorm[\cdot]_s$-norm for every
$s\newin\mn_0$. Then there exists a sequence of Banach spaces
$\{X_s\}_{s\in\mn_0}$, which satisfies (\ref{fx1})-(\ref{fx3}) and
such that $\{g_i|_{_{X_F}}\}_{i=1}^\infty$ is a pre-$F$-frame for
$X_F$ with respect to $\Theta_F$ with bounds $A_0, B_0,A_s=
B_s=1$, $s\newin\mn$.

If moreover, there exists an $F$-bounded projection from
$\Theta_F$ onto $M$, then $\{g_i|_{_{X_F}}\}_{i=1}^\infty$ is an
$F$-frame for $X_F$ with respect to $\Theta_F$.
\end{thm}

\bp For every $s\newin\mn$, define
\begin{equation}\label{normss}
X_s:=\left\{f\in X_0 \ : \ \{g_i(f)\}_{i=1}^\infty\in\Theta_s
\right\}, \; \|f \|_s:=\snorm[\{g_i(f)\}_{i=1}^\infty]_s.
\end{equation}
Let $s\in\mn$. Clearly, $X_s$ is a linear space, $\|\cdot\|_s$ is
a norm
 and furthermore,
$X_{s+1}\subseteq X_{s}$. Moreover, for every $f\in X_s$,
\begin{equation} \label{ns}
\|f \|_s=\snorm[\{g_i(f)\}_{i=1}^\infty]_{s}\geq
\snorm[\{g_i(f)\}_{i=1}^\infty]_{s-1}\left\{
\begin{array}{rl}
\geq A_0\|f\|_0, &  \hspace{-.05in} if \, s=1,\\
=\|f\|_{s-1}, & \hspace{-.05in} if \, s>1.
\end{array}
\right.\end{equation} Let $\{h_n\}_{n=1}^\infty$ be a Cauchy
sequence in $X_s$ and hence, by (\ref{ns}), a Cauchy sequence in
$X_0$. Therefore, $\{h_n\}_{n=1}^\infty$ converges to some element
$h\newin X_0$ and by the upper $\Theta_0$-frame inequality it
follows that
\begin{equation}\label{ghn}
\{g_i(h_n)\}_{i=1}^\infty \to \{g_i(h)\}_{i=1}^\infty\ \mbox{in}\
\Theta_0, \ \mbox{as} \ n\to\infty.
\end{equation}
By (\ref{normss}), $\{g_i(h_n)\}_{i=1}^\infty$, $n\newin\mn$,
 is a Cauchy sequence in $\Theta_s$
and it converges to some element $a$ in $\Theta_s$. Since
$\snorm[\cdot]_s\geq \,\snorm[\cdot]_0,$ the sequence
$\{g_i(h_n)\}_{i=1}^\infty$, $n\newin\mn$, converges to $a$ in
$\Theta_0$ and hence, by (\ref{ghn}),
 $\{g_i(h)\}_{i=1}^\infty=a\newin\Theta_s$. Therefore, $h\newin X_s$ and
$$\|h_n-h\|_s=
\snorm[\{g_i(h_n)\}_{i=1}^\infty-\{g_i(h)\}_{i=1}^\infty]_{s} \to
0,  \ \mbox{as} \ n\to\infty.$$ Thus, $X_s$ is a complete space.

By (\ref{normss}),  $\{g_i|_{{X_s}}\}_{i=1}^\infty$ is a tight
$\Theta_s$-frame for $X_s$ with frame bounds $A_s=B_s=1$.

Denote $X_F:=\sech[X]$. It remains to show the density of $X_F$ in
$X_s$, $s\newin\mn_0$. Since $\seqgr[g]$ is a $\Theta_0$-frame for
$X_0$, the operator $U_0:X_0\to\Theta_0$,
$U_0f=\{g_i(f)\}_{i=1}^\infty$, is injective. Observe that
$X_F=U_0^{-1}(M)$
 and hence $X_F\neq \{\nulel\}$,
because $M\neq \{\nulel\}$. Fix $s\newin\mn_0$. Let $f\newin X_s$
and $\varepsilon>0$. Since $M=U_0(X_F)$ is dense in $U_0(X_0)\cap
\Theta_s$, there exists $\tilde{f}\in X_F$ such that
$\snorm[U_0f-U_0\tilde{f}]_s<\varepsilon$ and hence,
$\|f-\tilde{f}\|_s\leq  \snorm[U_0(f-\tilde{f})]_s<\varepsilon$.

Therefore, $\{g_i|_{_{X_F}}\}_{i=1}^\infty$ is a pre-$F$-frame for
$X_F$ with respect to $\Theta_F.$

Assume now that there is an $F$-bounded projection from $\Theta_F$
onto $M$. Since $M=U_0(X_F)=\ru$, the pre-$F$-frame
$\{g_i|_{_{X_F}}\}_{i=1}^\infty$ is an $F$-frame.

\v Note that if $X_0 = X_s$ (as sets) for some $s \geq 1$ and thus
$X_0=X_t=X_s$ for any $0< t\leq s$, then $ \| \cdot \|_{X_0} $ and
$ \| \cdot \|_{X_t} $ are equivalent norms for any $0\leq t\leq
s$. This is a consequence of the Inverse Mapping Theorem. \ep

\v The following simple example shows a case when all the
assumptions of the above proposition are fulfilled.

\begin{exmp}
{\rm
 Let $\{\Theta_s\}_{ s\in\mn_0}$ be a sequence of $CB$-spaces
 which satisfies (\ref{fx1})-(\ref{fx3}) and let $X_0=\Theta_0$.
 Let $\{g_i\}_{i=1}^\infty\newin
(\Theta_0^*)^\mn$ be the sequence of the coordinate functionals,
associated to the canonical basis $\seqgr[e]$ of \,$\Theta_0$,
i.e. $g_j(\seqgr[x])=x_j$, $j\newin\mn$. Clearly,
$\{g_i\}_{i=1}^\infty$ is a $\Theta_0$-frame for $\Theta_0$ with
bounds $A_0=B_0=1$. For every $s\in\mn$ the space $X_s$,
constructed in the proof of Theorem \ref{constructx}, coincides
with $\Theta_s$ and $\cap_{s\in} X_s=\Theta_F \neq\{\nulel\}$. In
this case $M= \Theta_F=\ru$ and $\{g_i|_{_{X_F}}\}_{i=1}^\infty$
is an $F$-frame for $X_F$ with respect to $\Theta_F$. }
\end{exmp}

\subsection{Construction of $\{\Theta_s\}_{s\in\mn_0}$ - restrictive case}

In the subsections which are to follow we consider the problem
opposite to the one in the previous subsection: given
$\{X_s\}_{s\in\mn_0},$ construct $\{\Theta_s\}_{s\in\mn_0}$.

\begin{thm}\label{constrth1}
Let $\{X_s\}_{s\in\mn_0}$ be a sequence of Banach spaces, which
satisfies (\ref{fx1})-(\ref{fx3}). Let $\Theta$ be a $BK$-space
and let $\seqgr[g]\newin (X_0^*)^\mn$ satisfy the lower
$\Theta$-frame condition for $X_0$, i.e. there exists a constant
$A>0$ such that $\{g_i(f)\}_{i=1}^\infty \in\Theta$ and
$A\|f\|_0\leq  \snorm[\{g_i(f)\}_{i=1}^\infty]$ for every $f\newin
X_0$. Then the following holds.

{\rm (a)} There exists a sequence $\{\Theta_s\}_{s\in{\mn_0}}$ of
$BK$-spaces, which satisfies (\ref{fx1})-(\ref{fx3}), such that
$\{g_i|_{_{X_F}}\}_{i=1}^\infty$ is an $F$-frame for $X_F$ with
respect to $\Theta_F$ with bounds $A_s=B_s=1$, $s\newin\mn_0$.

{\rm (b)} If $\seqgr[g]$ has a biorthogonal sequence
$\seqgr[f]\newin (X_F)^\mn$, then the spaces $\Theta_s,$
$s\newin\mn_0$,
 constructed in (a),  are $CB$-spaces if and only if
$\seqgr[f]$ is a Schauder basis of $X_s$, $s\newin\mn_0$.
\end{thm}
\bp (a) The lower $\Theta$-frame inequality implies the unique
correspondence of $f$ to $\{g_i(f)\}_{i=1}^\infty$. This leads to
the following definition:
\begin{equation}\label{ths}
\Theta_s:=\left\{ \{g_i(f)\}_{i=1}^\infty \, : \, f\in
X_s\right\}\subset \Theta,\;\;
\snorm[\{g_i(f)\}_{i=1}^\infty]_{s}:=\|f\|_s, \, s\in\mn_0.
\end{equation}
Let $s\newin\mn_0$. Clearly, $\Theta_s$ is a linear space,
$\snorm[\cdot]_{s}$ is a norm in $\Theta_s$ and
$\Theta_{s+1}\subseteq \Theta_{s}$. It is easy to see that
$\Theta_s$ is complete and $\snorm[\cdot]_{s+1}\geq
\snorm[\cdot]_{s}$.

For $k\newin\mn$, let $E_k$ denote the $k$-th coordinate
functional $E_k(\{c_i\}_{i=1}^\infty)=c_k$. For every $f\newin
X_s$ one has
\begin{eqnarray*}
\left|E_k\left(\{g_i(f)\}_{i=1}^\infty\right)\right|&=&
|g_k(f)|\leq \|g_k\|_{X_0^*}\, \|f\|_0\leq  \|g_k\|_{X_0^*}\, \|f\|_s\\
&=& \|g_k\|_{X_0^*}\, \snorm[\{g_i(f)\}_{i=1}^\infty]_s,
\end{eqnarray*}
which implies that $E_k$ is continuous on $\Theta_s$,
$s\newin\mn_0$. Thus, $\Theta_s$ is a $BK$--space. Moreover,
$\Theta_F=\cap_{s\in\mn_0}\Theta_s=\{\{g_i(f)\}_{i=1}^\infty :
f\newin X_F\}\neq\{\nulel\}$, because $X_F\neq\{\nulel\}$ and
$g_i\neq 0$ for at least one $i\in \mn$.
 Furthermore, for every $s\newin\mn_0$, $\Theta_F$ is dense in $\Theta_s$, because $X_F$ is dense in $X_s$, which is isomorphic to $\Theta_s$.

By (\ref{ths}), $\{g_i|_{_{X_F}}\}_{i=1}^\infty$ is a
pre-$F$-frame for $X_F$ with respect to $\Theta_F$ with bounds
$A_s=B_s=1$, $s\newin\mn_0$. Since
$\Theta_F=\{\{g_i(f)\}_{i=1}^\infty : f\newin X_F\}$, (\ref{ths})
implies that the operator $V: \Theta_F\to X_F$ given by
$V(\{g_i(f)\})=f$ is $s$-bounded for every $s\newin\mn_0$. Thus,
$\{g_i|_{_{X_F}}\}_{i=1}^\infty$ is an $F$-frame for $X_F$ with
respect to $\Theta_F$.

(b) Let $\seqgr[g]$ have a biorthogonal sequence $\seqgr[f]\newin
(X_F)^\mn$. Consider $\{\Theta_s\}_{s\in\mn_0}$ constructed in
(a). Fix $s\newin \mn_0$ and $f\newin X_s$. By the definition of
$\Theta_s$, for every $k\newin\mn_0$, the $k$-th canonical vector
$e_k=\{g_i(f_k)\}_{i=1}^\infty$ belongs to $\Theta_s$. By the
biorthogonality, we have $$ \sum_{i=1}^n
g_i(f)e_i=\{g_i(g_1(f)f_1+g_2(f)f_2+\ldots+g_n(f)f_n
)\}_{i=1}^\infty, \ n\in\mn,$$ and hence
\begin{eqnarray*}
\snormb[\{g_i(f)\}_{i=1}^\infty-\sum_{i=1}^n g_i(f)e_i]_s
&=&\snormb[\left\{g_i\left(f - \sum_{i=1}^n
g_i(f)f_i\right)\right\}_{i=1}^\infty]_s\\
&=& \left\|f - \sum_{i=1}^n g_i(f)f_i \right\|_s.
\end{eqnarray*}
This implies that
$$\sum_{i=1}^n g_i(f)e_i\to
\{g_i(f)\}_{i=1}^\infty \ \mbox{as}\ n\to\infty \ \mbox{if and
only if}\ \sum_{i=1}^n g_i(f)f_i \to f \ \mbox{as}\ n\to\infty.$$

 The conclusion now follows from the fact that the sequence
$\seqgr[f]$, which has $\seqgr[g]$ as a biorthogonal sequence, is
a Schauder basis of $X_s$ if and only if $f=\sumii g_i(f) f_i$ for
every $f\newin X_s$ (see e.g. \cite{Singer}, also \cite{Heil}).
\ep

The example which is to follow illustrates the above theorem.

\begin{exmp}\label{herm}
{\rm
Let $\{h_n\}_{n=0}^\infty$ be the Hermite basis of $L^2(\mr)$:
$$h_n(x)=\frac{1}{(2^n\, n! \sqrt\pi)^{1/2}} e^{-\frac{x^2}{2}} H_n(x), \ \
x\in\mr,
$$
where $H_n(x)=(-1)^n \, e^{x^2} (d/dx)^n \, (e^{-x^2})$,
$x\newin\mr$, is the $n$-th Hermite polynomial, $n\newin\mn_0$.
The harmonic oscilator $\mathcal{R}:=-\frac{d^2}{dx^2}+x^2$ has
eigenvalues $\lambda_n=2n+1$, $n\newin\mn_0$. Define
$$X_s:=\{f\in L^2(\mr)\ :  \ \mathcal{R}^sf\in L^2\},\ \ \|f\|_s:=\|\mathcal{R}^sf\|_{L^2}, \ s\in \mn_0. $$
Then $\{X_s\}_{s\in\mn_0}$ is a sequence of Hilbert spaces, which
satisfies (\ref{fx1})-(\ref{fx3}) and
$\{\frac{1}{(2n+1)^s}h_n\}_{n=0}^\infty$ is an orthonormal basis
for $X_s$ for every $s\in\mn_0$. Moreover,
$\{\Theta_s\}_{s\in{\mn_0}}$, defined by (\ref{ths}), is a
sequence of $CB$-spaces, which satisfies (\ref{fx1})-(\ref{fx3})
and $\{g_i|_{_{X_F}}\}_{i=1}^\infty$ is an $F$-frame for $X_F$
with respect to $\Theta_F$.
}
\end{exmp}

As it was already noted, the $CB$-property is crucial for the
series expansions via Banach frames. Because of that, in Theorem
\ref{constrth1}\,(b) we point out construction of $CB$-spaces
$\Theta_s$, $s\newin\mn_0$.
 However, in order to have
the canonical vectors in $\Theta_s$ (which contains only sequences
of the form $\{g_i(f)\}_{i=1}^\infty$), we assume the existence of
a biorthogonal sequence $\seqgr[f]$, which restricts us to the
cases when $\seqgr[f]$ is a Schauder basis for $X_s$. Note that
when $X_s$ is reflexive, $\seqgr[f]$ is a Schauder basis of $X_s$
if and only if its biorthogonal sequence $\seqgr[g]$ is a Schauder
basis of $X_s^*$ \cite[Corollary 8.2]{Heil}. This restricts
 us to the use of sequences $\seqgr[g]$,
 which
 are Schauder basis. In order to avoid the use of a biorthogonal sequence, in the
following subsection we construct larger sequence spaces. Thus,
our forthcoming Theorems \ref{constructone} and
\ref{constructthetas} are related to the cases when $\seqgr[g]$ is
not necessarily a Schauder basis.

\subsection{Construction of $\{\Theta_s\}_{s\in\mn_0}$ - more general case} \label{subs43}

In the sequel we will use the following notation related to a
sequence $c=\seqgr[c]$:
$$c^{(n)}:=\{\underbrace{0,\ldots,0}_n,c_{n+1},c_{n+2},c_{n+3},\ldots\}=c-\sum_{i=1}^n c_i e_i, \ n\in\mn,$$
 \hspace{.55in} $c^{(n)}_i$:= the $i$-th coordinate of $c^{(n)}$,
$i\in\mn$.

\v For the construction in Theorem \ref{constructthetas} we need
the following theorem.

\begin{thm} \label{constructone}
Let $\Theta\neq\{\nulel\}$ be a solid $BK$-space,
$X\neq\{\nulel\}$ be a reflexive Banach space and $\seqgr[g]\newin
(X^*)^\mn$ be a $\Theta$-Bessel sequence for $X$ with  bound
$B\leq 1$ such that $0<\|g_i\|\leq 1$, $i\newin\mn$. For every
$c=\seqgr[c]\newin \Theta$, denote
\begin{equation} \label{setm}
M^c:=\{f\in X \ : \ |c_i|\leq |g_i(f)|, \ i\in \mn \}
\end{equation} and define
\begin{equation} \label{setm2}
 \widetilde{\Theta} := \left\{ c\in\Theta \ : \ M^{c}\neq
\varnothing\right\}, \ \, \snorm[c]_{\widetilde{\Theta}}
:=\inf\left\{ \|f\| : f\in M^{c}\right\}.
\end{equation}
Assume

\vspace{.1in} \noindent $(\mathcal{A}_1): \ (\forall
\,c\in\widetilde{\Theta})\ (\forall \, d\in \widetilde{\Theta}) \
(\forall\, f\in M^c) \ (\forall\, h\in M^d) \Rightarrow $

\vspace{.1in} \hspace{1.7in} $ (\exists \ r\in M^{c+d})\ \
(\|r\|\leq \|f\|+\|h\|).$

\noindent Then: \begin{itemize}
\item[{\rm(a)}]
$\widetilde{\Theta}\,(\subseteq \Theta)$ is a solid $BK$-space
with
$\snorm[\cdot]_{\Theta}\leq\snorm[\cdot]_{\widetilde{\Theta}}$ and
$\{g_i\}_{i=1}^\infty$ is a $\widetilde{\Theta}$-Bessel sequence
for $X$ with bound $\widetilde{B}=1$.
\item[{\rm(b)}] If $\seqgr[g]$ is a $\Theta$-frame for
$X$, then $\{g_i\}_{i=1}^\infty$ is a $\widetilde{\Theta}$-frame
for $X$.
\item[{\rm(c)}] If $\seqgr[g]$ is a Banach frame for
$X$ with respect to $\Theta$, then $\{g_i\}_{i=1}^\infty$ is a
Banach frame for $X$ with respect to $\widetilde{\Theta}$.
\item[{\rm(d)}] $\widetilde{\Theta}$ is a $CB$-space if and only if

\vn $(\mathcal{A}_2): \ (\forall c\in \widetilde{\Theta}) \
(\forall\, \varepsilon>0) \ (\exists \, k\in\mn) \  (\exists
\,f\in M^{c^{(k)}})  \ ( \|f\|<\varepsilon).$
\end{itemize}
\end{thm}

\bp (a) First note that $\widetilde{\Theta}\neq\{\nulel\}$. For
any $i\newin \mn$, $g_i$ is not the null functional and hence
there exists $f\newin X$ such that $1\leq |g_i(f)|$. This implies
that the $i$-th canonical vector $e_i$ belongs to
$\widetilde{\Theta}$.

Let $c=\seqgr[c]\in \widetilde{\Theta}$ and $f\newin M^{c}$. Since
$\Theta$ is solid, (\ref{setm}) implies that
$$
\snorm[\{c_i\}_{i=1}^\infty]_\Theta\leq
\snorm[\{g_i(f)\}_{i=1}^\infty]_\Theta\leq B \,\|f\|\leq \|f\|.
$$
Therefore,
\begin{equation}\label{normnerav}
\snorm[\{c_i\}_{i=1}^\infty]_\Theta \leq
 \snorm[\{c_i\}_{i=1}^\infty]_{\widetilde{\Theta}}.\end{equation}

Let $c\in\widetilde{\Theta}$ and $d\in\widetilde{\Theta}$. By
$(\mathcal{A}_1)$, the sequence $c+d$ belongs to
$\widetilde{\Theta}$; clearly $\lambda c\in\widetilde{\Theta}$,
$\lambda\in\mr$. Thus $\widetilde{\Theta}$ is a linear space. In
order to prove that $\snorm[\cdot]_{\widetilde{\Theta}}$ is a
norm, let us first note that
$\snorm[\{c\}_{i=1}^\infty]_{\widetilde{\Theta}}=0\
\mbox{implies}\ c_i=0, i\newin\mn.$ This is a consequence of the
inequality (\ref{normnerav})
 and the fact that $\snorm[\cdot]_\Theta$ is a norm. The equality
$\snorm[\lambda c]_{\widetilde{\Theta}}=|\lambda|\,
\snorm[c]_{\widetilde{\Theta}}$ follows from the fact that
$$M^{\lambda c}=\{ \lambda f \ : \ f\in M^{c}\}, \ \lambda\neq 0.$$
For the triangle inequality, fix arbitrary
$c\newin\widetilde{\Theta}$, $d\newin\widetilde{\Theta}$. For
every $f\newin M^c$ and every $h\newin M^d$, by ($\mathcal{A}_1$)
we can choose an element $r^{f,h}\newin M^{c+d}$ such that
$\|r^{f,h}\|\leq \|f\|+\|h\|$. Then
\begin{eqnarray*}
\snorm[c+d]_{\widetilde{\Theta}}&=& \inf\left\{ \|r\| \ :\ r\in
M^{c+d}\right\} \leq \inf\left\{ \|r^{f,h}\| \ : \
f \in M^c, h\in M^d \right\}\\
&\leq& \inf\left\{ \|f\| \ :\ f\in M^{c}\right\} + \inf\left\{
\|h\| \ :\
h\in M^{d} \right\}\\
&=&\snorm[c]_{\widetilde{\Theta}}+\snorm[d]_{\widetilde{\Theta}}.\end{eqnarray*}

It remains to show that $\widetilde{\Theta}$ is complete. Let
$c^\nu=\{c_i^\nu\}_{i=1}^\infty, \nu\in\mn,$ be a Cauchy sequence
in $\widetilde{\Theta}$. Fix $\varepsilon
>0$. There exists $\nu_0(\varepsilon)$ such that for
every $\mu,\nu\in\mn$, $\mu\geq\nu_0$, $\nu\geq \nu_0$, there
exists $f^{\mu,\nu}\in M^{c^\mu - c^\nu}$, such that
$$
 \|f^{\mu,\nu}\|<\varepsilon \ \mbox{and} \ |c_i^\mu -
c_i^\nu|\leq  |g_i(f^{\mu,\nu})|, \
 i\in\mn.
$$
By (\ref{normnerav}), $c^\nu, \nu\newin\mn$, is a Cauchy sequence
in $\Theta$ and hence there exists $c=\seqgr[c]\in\Theta$ such
that
$$\snorm[\{c_i^\nu-c_i\}_{i=1}^{\infty}]_\Theta\to 0 \ \mbox{as} \ \nu\to \infty.$$
 Let $i\in\mn$. Since $\Theta$ is $BK$-space, it follows that
 $c_i^\nu\to c_i$ as $\nu\to\infty$.
Fix $\nu\geq \nu_0$. We will find an element which belongs to
$M^{c-c^\nu}$. Since $\|f^{\mu,\nu}\|<\varepsilon $ for every
$\mu\geq\nu_0$, the scalar sequence
$\{\|f^{\mu,\nu}\|\}_{\mu=\nu_0}^\infty$ is bounded and hence it
contains a convergent subsequence
$\{\|f^{\mu_k,\nu}\|\}_{k=1}^\infty$; denote its limit by $a^\nu$.
Since $X$ is reflexive and the sequence
$\{f^{\mu_k,\nu}\}_{k=1}^\infty$ is norm-bounded, by
\cite[Corollary 1.6.4]{AK} there exists a subsequence
$\{f^{\mu_{k_n},\nu}\}_{n=1}^\infty$ which converges weakly to
some element $F^\nu\newin X$. Therefore,
$$\|F^\nu\|\leq
\lim \inf
\|f^{\mu_{k_n},\nu}\|=\lim_{n\to\infty}\|f^{\mu_{k_n},\nu}\|=a^\nu.$$
Since for every $n\in\mn$, $\|f^{\mu_{k_n},\nu}\|<\varepsilon$, it
follows
$$a^\nu = \lim_{n\to\infty} \|f^{\mu_{k_n},\nu}\|\leq
\varepsilon \ \mbox{and hence} \ \|F^\nu\|\leq\varepsilon.$$ Fix
$i\newin\mn$. For every $n\newin\mn$, $$ |c_i^{\mu_{k_n}} -
c_i^\nu|\leq |g_i(f^{\mu_{k_n},\nu})|$$ and taking limit as
$n\to\infty$, we obtain
$$ |c_i - c_i^\nu|\leq
|g_i(F^\nu)|.$$ Thus, $F^\nu$ belongs to $M^{c-c^\nu}$. This
implies that $c-c^\nu\in \widetilde{\Theta}$ and therefore,
$c\in\widetilde{\Theta}$. Moreover,
$$\snorm[\{c_i-c_i^\nu\}_{i=1}^\infty]_{\widetilde{\Theta}}=\inf \{\|f\| \ : \ f\in
M^{c-c^\nu}\} \leq  \|F^\nu\|\leq \varepsilon, \ \nu\geq \nu_0.$$
This concludes the proof that $\widetilde{\Theta}$ is complete.

Since the coordinate functionals are continuous on $\Theta$ and
$\snorm[\cdot]_\Theta\leq\snorm[\cdot]_{\widetilde{\Theta}}$, they
are continuous on $\widetilde{\Theta}\subseteq\Theta$ and thus
$\widetilde{\Theta}$ is a $BK$-space. Moreover, the space
$\widetilde{\Theta}$ is solid. Indeed, let $c=\seqgr[c]\in
\widetilde{\Theta}$ and $d=\seqgr[d]$ be such that $|d_i|\leq
|c_i|$, $i\newin\mn$. Since $\Theta$ is solid, $\seqgr[d]$ belongs
to $\Theta$. Moreover, $M^{d}\supseteq M^{c}\neq \varnothing$.
Thus, $\seqgr[d]\in\widetilde{\Theta}$ and
$\snorm[d]_{\widetilde{\Theta}}\leq\snorm[c]_{\widetilde{\Theta}}$.

For every $f\newin X$, the set $M^{\{g_i(f)\}_{i=1}^\infty}$
contains $f$ and thus
$\{g_i(f)\}_{i=1}^\infty\newin\widetilde{\Theta}$; moreover,
$\snorm[\{g_i(f)\}_{i=1}^\infty]_{\widetilde{\Theta}}\leq\|f\|$.
Therefore, $\{g_i\}_{i=1}^\infty$ is a $\widetilde{\Theta}$-Bessel
sequence for $X$ with bound $\widetilde{B}=1$.

(b) Let $A$ be a lower $\Theta$-frame bound for the $\Theta$-frame
$\seqgr[g]$. Then
$$A\|f\|\leq \snorm[\{g_i(f)\}_{i=1}^\infty]_{\Theta}\leq\snorm[\{g_i(f)\}_{i=1}^\infty]_{\widetilde{\Theta}}, \ f\in X,$$
which implies that $\seqgr[g]$ satisfies the lower
$\widetilde{\Theta}$-frame inequality.

(c) Let $V:\Theta\to X$ be a bounded operator such that
$V(\{g_i(f)\}_{i=1}^\infty)=f$, $f\newin X$. Consider
$\widetilde{V}:=V\vert_{\widetilde{\Theta}}$. For every $c\in
{\widetilde{\Theta}}$,
$$\|\widetilde{V} c\|=\|V c\|\leq \|V\|\,\snorm[c]_\Theta \leq \|V\|\,
\snorm[c]_{\widetilde{\Theta}},
 $$
which implies that $\widetilde{V}$ is bounded on
$\widetilde{\Theta}$. Therefore, $\seqgr[g]$ is a Banach frame for
$X$ with respect to $\widetilde{\Theta}$.

 (d) That all the canonical vectors belong to
$\widetilde{\Theta}$, is shown in (a). For
$c\newin\widetilde{\Theta}$, observe that $M^{c^{(n)}}\subseteq
M^{c^{(n+1)}}$, $n\in\mn$, and hence
$\snorm[c^{(n+1)}]_{\widetilde{\Theta}}\leq
\snorm[c^{(n)}]_{\widetilde{\Theta}}$, $n\in\mn$.

Assume now that $(\mathcal{A}_2)$ holds. Let
$c\in\widetilde{\Theta}$ and $\varepsilon>0$. By
$(\mathcal{A}_2)$, there exists $k\in\mn$ such that
$\snorm[c^{(k)}]_{\widetilde{\Theta}}<\varepsilon$ and hence
$\snorm[c^{(n)}]_{\widetilde{\Theta}}<\varepsilon$ for every
$n\geq k$. Therefore, $c^{(n)}\to 0$ in ${\widetilde{\Theta}}$ as
$n\to\infty$, which implies that $\sum_{i=1}^n c_i e_i\to c$ in
$\widetilde{\Theta}$ as $n\to\infty.$

Assume now that ${\widetilde{\Theta}}$ is a $CB$-space and let
$c\in \widetilde{\Theta}$. For every $\varepsilon >0$ there exists
$k\in\mn$ such that
$\snorm[c^{(k)}]_{\widetilde{\Theta}}<\varepsilon$. This implies
that there exists $f\in M^{c^{(k)}}$ with $\|f\|<\varepsilon$.
 \ep

\v In order to show that conditions $(\mathcal{A}_1)$ and
$(\mathcal{A}_2)$ are not artificial ones, we will prove their
validity in a simple case when a frame is obtained from an
orthonormal basis of a Hilbert space.

\begin{prop}\label{primer1}
Let $(X, \<\cdot,\cdot\>)$ be a Hilbert space, $\seqgr[e]$ be an
orthonormal basis for $X$ and $\Theta=\ell^2$. Let
$\seqgr[g]\newin (X^*)^\mn$ be defined by
$$g_1(f):=\<f, e_1\> \ \mbox{and}\ \, g_i(f):=\<f, e_{i-1}\>, \ i=2,3,4,\ldots.$$
Clearly, $\seqgr[g]$ is a Banach frame for $X$ with respect to
$\ell^2$, which is not a Schauder basis of $X$. Let
$\widetilde{\Theta}$ be defined by (\ref{setm2}).
  Then $\widetilde{\Theta}$ is a $CB$-space and $\seqgr[g]$ is a Banach frame for $X$ with respect to $\widetilde{\Theta}$.
 \end{prop}
 \bp
Let $c=\seqgr[c]\in\widetilde{\Theta},
d=\seqgr[d]\in\widetilde{\Theta}$. Denote
$$\widetilde{c}:=\max(|c_1|, |c_2|), \ \, \widetilde{d}:=\max(|d_1|,
|d_2|).$$ Fix $f\in M^c$ and $h\in M^d$. By (\ref{setm}),
\begin{equation} \label{nova}
|c_i|\leq |g_i(f)|,\   |d_i|\leq |g_i(h)|, \ i\in \mn.
\end{equation}
 Let us find $r^{f,h}\in M^{c+d}$ such that $\|r^{f,h}\|\leq \|f\|+\|h\|$.
Consider
$$r^{f,h}:= m e_1 + |c_3+d_3| e_2 + |c_4+d_4| e_3 + |c_5+d_5| e_4
+\ldots \in X,$$ where $m:= \max(|c_1+d_1|, |c_2+d_2|).$ It is
clear that
$$|c_i+d_i|=\<r^{f,h},e_{i-1}\>=g_i(r^{f,h}), \ i\geq 3,$$
$$|c_i+d_i|\leq
m = \<r^{f,h},e_1\>=g_i(r^{f,h}), \ i=1,2,$$ which implies that
$r^{f,h}\newin M^{c+d}$. By (\ref{nova}), $$\widetilde{c}\leq
|\<f,e_1\>|, \ \ |c_i|\leq |\<f,e_{i-1}\>|, \ i\geq 3,$$
$$ \widetilde{d}\leq |\<h,e_1\>|, \ \ |d_i|\leq |\<h,e_{i-1}\>|,\ i\geq 3.$$
Using the fact that $m\leq \widetilde{c}+\widetilde{d}$,
$|c_i+d_i|\leq |c_i|+|d_i|$, and the solidity of $\ell^2$, we
obtain that
\begin{eqnarray*}
\|r^{f,h}\| &\leq & \|\{\widetilde{c}+\widetilde{d}, |c_3|+|d_3|,
|c_4|+|d_4|,
\ldots \} \|_{\ell^2}\\
& \leq & \|\{\widetilde{c}, |c_3|, |c_4|, \ldots \}\|_{\ell^2}
+\|\{\widetilde{d},|d_3|,|d_4|, \ldots \}
\|_{\ell^2}\\
& \leq & \|\{|\<f,e_i\>|\}_{i=1}^\infty\|_{\ell^2}
+\|\{|\<h,e_i\>| \}_{i=1}^\infty \|_{\ell^2}= \|f\|+\|h\|.
\end{eqnarray*}
Therefore, $(\mathcal{A}_1)$ is fulfilled.

Consider now $c=\seqgr[c]\in\widetilde{\Theta}$. Fix
$\varepsilon>0$ and find $k\newin\mn$, $k>1,$ such that
$\sum_{i=k+1}^\infty |c_i|^2<\varepsilon$. Let
$$b:=\{\underbrace{0, \ldots, 0}_{k-1}, |c_{k+1}|,|c_{k+2}|,|c_{k+3}|,\ldots \}.$$
There exists $h\newin X$ such that $b_i=\<h,e_i\>, i\in\mn, \ \,
\mbox{and} \ \, \|h\|^2=\sum_{i=1}^\infty |b_i|^2.$ Recall,
$$c^{(k)}=\{\underbrace{0,\ldots,0}_k,c_{k+1},c_{k+2},c_{k+3},\ldots\}.$$
Thus we have
$$|c^{(k)}_1|=0\leq |g_{1}(h)| \ \, \mbox{and} \ \,
|c^{(k)}_i|=b_{i-1}=\<h,e_{i-1}\>= g_{i}(h), \ \ i\geq 2,$$ which
implies that $h\in M^{c^{(k)}}$. Moreover,
$\|h\|^2=\sum_{i=k+1}^\infty |c_i|^2<\varepsilon $ and hence
$(\mathcal{A}_2)$ is fulfilled.

The conclusion now follows from Theorem \ref{constructone}. \ep

\begin{thm}\label{constructthetas}
  Let $\Theta\neq\{\nulel\}$ be a solid
$BK$-space and $\{X_s\}_{s\in\mn_0}$ be a sequence of reflexive
Banach spaces which satisfies (\ref{fx1})-(\ref{fx3}). Let
$\seqgr[g]\newin (X_0^*)^\mn$ be a $\Theta$-Bessel sequence for
$X_0$ with bound $B\leq 1$ such that $0<\|g_i\|\leq 1$,
$i\newin\mn$. For every $s\newin\mn_0$ and every
$c=\seqgr[c]\newin \Theta$, denote
\begin{equation}\label{mx2}
M^{c}_s:=\{f\in X_s \ : \ |c_i|\leq |g_i(f)|, \ i\in \mn \}
\end{equation}
and define
 \begin{equation}\label{ts2}
 {\Theta}_s := \left\{ c\in\Theta \ : \ M^c_s\neq
\varnothing\right\}, \ \, \snorm[c]_{s} :=\inf\left\{ \|f\|_s :
f\in M^c_s\right\}.
 \end{equation}
Consider the following conditions:

\vspace{.1in}\noindent $(\mathcal{A}_1^\prime):$ \ \ \ \ \
 $(\forall c\in\Theta_s)\ \ (\forall d\in \Theta_s)\ \ (\forall f\in
M^c_s) \ \ (\forall h\in M^d_s) \ \Rightarrow $

\vspace{.1in} \centerline{$(\exists \,r\in M^{c+d}_s)\ \
(\|r\|_s\leq \|f\|_s+\|h\|_s).$}

\vspace{.1in}\noindent $(\mathcal{A}_2^\prime): \
 \ (\forall c\newin \Theta_s) \  (\forall\,
\varepsilon>0)\ \ (\exists \, k\newin\mn) \  \, (\exists \,f\newin
M^{c^{(k)}}_s) \, \ (\|f\|_s<\varepsilon).$

\vspace{.1in}\noindent $(\mathcal{A}_3^\prime):$ \ There exists
$A_s\in(0,1]$ such that for every $f\in X_s$ one has
\begin{equation}\label{fj1}
\widetilde{f}\in M_s^{\{g_i(f)\}_{i=1}^\infty} \Rightarrow
A_s\|f\|_s\leq \|\widetilde{f}\|_s.
\end{equation}

Assume that $(\mathcal{A}_1^\prime)$ holds for every $s\in\mn_0$.
Then:

\begin{itemize}
\item[${\rm(a)}$] $\{\Theta_s\}_{s\in\mn_0}$ is a sequence of solid
$BK$-spaces with the properties (\ref{fx1})-(\ref{fx2}) such that
$\{g_i|_{X_s}\}_{i=1}^\infty$ is a $\Theta_s$-Bessel sequence for
$X_s$ with bound $B_s=1$, $s\newin\mn_0$.
\item[${\rm(b)}$] For any $s\newin\mn$,
$\Theta_s$ is a $CB$-space if and only if $(\mathcal{A}_2^\prime)$
holds.
\item[${\rm(c)}$] For any $s\newin\mn$,
 $\{g_i|_{X_s}\}_{i=1}^\infty$ is a
$\Theta_s$-frame for $X_s$ if and only if $(\mathcal{A}_3^\prime)$
holds.

If $(\mathcal{A}_3^\prime)$ holds with $A_s=1$, then
$\{g_i\vert_{X_s}\}$ is a tight $\Theta_s$-frame for $X_s$.
\end{itemize}

\end{thm}

\bp $(a)$ Let $s\in\mn_0$. Since $\seqgr[g]\newin (X_0^*)^\mn$ and
$\seqgr[g]$ is a $\Theta$-Bessel sequence for $X_0$ with bound
$B$, (\ref{fx2}) implies that
$\{g_i\vert_{X_s}\}_{i=1}^\infty\newin (X_s^*)^\mn$ and
$\{g_i\vert_{X_s}\}_{i=1}^\infty$ is a $\Theta$-Bessel sequence
for $X_s$ with bound $B$. Now Proposition \ref{constructone},
applied to $X=X_s$ and $\{g_i\vert_{X_s}\}_{i=1}^\infty$, implies
that $\Theta_s$ is a solid $BK$-space and
$\{g_i\vert_{X_s}\}_{i=1}^\infty$ is a $\Theta_s$-Bessel sequence
for $X_s$ with bound $B_s=1$.

 It is clear that $\Theta_{s+1}\subseteq \Theta_s$ and
$\snorm[c]_s\leq\snorm[c]_{s+1}$, $ c\newin \Theta_{s+1},
s\newin\mn$, because $M^{c}_{s+1}\subseteq M^{c}_{s}$ and
$\|\cdot\|_s\leq \|\cdot\|_{s+1}$.

For every $f\newin X_s$, the set $M^{\{g_i(f)\}_{i=1}^\infty}_s$
contains  $f$ and thus
$$\Theta_s\supseteq \left\{ \{g_i(f)\} \, : \, f\in X_s\right\}.$$
Therefore, $$\Theta_F:=\cap_{s\in\mn_0}\Theta_s \supseteq
\{\{g_i(f)\} : f \in X_F\}\neq \{\nulel\},$$ because $X_F\neq
\{\nulel\}$ and $g_i\neq \nulel$, $i\newin\mn$.

$(b)$ Validity of the statement follows from Proposition
\ref{constructone}, applied to $X=X_s$ and to the $\Theta$-Bessel
sequence $\{g_i\vert_{X_s}\}_{i=1}^\infty$ for $X_s$.

$(c)$ Let $s\in\mn_0$. If $(\mathcal{A}_3^\prime)$ holds, then
$$A_s\|f\|_s\leq \inf\{\|\widetilde{f}\|_s \ : \ \widetilde{f}\in
M_s^{\{g_i(f)\}_{i=1}^\infty}\}=\snorm[\{g_i(f)\}]_{s}, \ f\in
X_s.$$ Conversely, if $\{g_i\vert_{X_s}\}_{i=1}^\infty$ is a
$\Theta_s$-frame for $X_s$ with a lower bound $A_s\in (0,1]$, then
for every $f\in X_s$ one has
$$A_s\|f\|_s
\leq \snorm[\{g_i(f)\}]_{s} \leq \|\widetilde{f}\|, \
\widetilde{f}\in M_s^{\{g_i(f)\}_{i=1}^\infty}.$$

Additionally, if $(\mathcal{A}_3^\prime)$ holds with $A_s=1$, then
$B_s=1$ implies that $\{g_i\vert_{X_s}\}_{i=1}^\infty$ is a tight
$\Theta_s$-frame for $X_s$. \ep

 \v Direct consequences of Theorem \ref{constructthetas} are given in
 the next corollary.
\begin{cor}\label{c1}
\noindent {\rm(a)} If $(\mathcal{A}_1^\prime)$ and
$(\mathcal{A}_2^\prime)$ are satisfied for every $s\in\mn_0$, then
$\{\Theta_s\}_{s\in\mn_0}$ is a sequence of solid $CB$-spaces with
the properties (\ref{fx1})-(\ref{fx3}) such that
$\{g_i|_{X_F}\}_{i=1}^\infty$ is an $F$-Bessel sequence for $X_F$
with respect to $\Theta_F$.

\vspace{.1in}\noindent {\rm(b)} If $(\mathcal{A}_1^\prime)$,
$(\mathcal{A}_2^\prime)$  and $(\mathcal{A}_3^\prime)$ are
satisfied for every $s\in\mn_0$, then
$\{g_i|_{X_F}\}_{i=1}^\infty$ is a pre-$F$-frame for $X_F$ with
respect to $\Theta_F$.
\end{cor}

Similarly as in Proposition \ref{primer1}, we show in Proposition
\ref{hilbcor} that conditions $(\mathcal{A}_1^\prime)
-(\mathcal{A}_3^\prime)$ are not artificial. For an illustration
of the assumptions of Proposition \ref{hilbcor}, see Example
\ref{herm} with $e_i=h_{i-1}$ and $a_{i,s}={(2i-1)^s}$,
$i\newin\mn$, $s\newin\mn$.

\begin{prop}\label{hilbcor}
Let $(X_0, \<\cdot,\cdot\>_0)$ be a Hilbert space and let
$\seqgr[e]$ denote an orthonormal basis for $X_0$. For given
number sequences $\{a_{i,s}\}_{i=1}^\infty$, $s\newin\mn$, with
$1\leq a_{i,s}\leq a_{i,s+1}$, $i\newin\mn$, $s\newin\mn$, define
$$
X_s:=\left\{f\in X_0 \ : \ \{a_{i,s} \<f,e_i\>_0\}_{i=1}^\infty
\in \ell^2\right\}, \ \ \<f,h\>_s:=\sumii a_{i,s}^2 \<f,e_i\>_0
\,\<e_i, h\>_0.
$$
Let $\Theta=\ell^2$ and $\seqgr[g]\newin (X_0^*)^\mn$ be defined
by
$$g_1(f):=\<f, e_1\>_0 \ \mbox{and}\ \, g_i(f):=\<f, e_{i-1}\>_0, \ i=2,3,4,\ldots.$$
Then $\{X_s\}_{s\in\mn_0}$ is a sequence of Hilbert spaces, which
satisfies (\ref{fx1})-(\ref{fx3}); $\{\Theta_s\}_{s\in\mn_0}$,
constructed by (\ref{mx2}) and (\ref{ts2}), is a sequence of
$CB$-spaces, which satisfies (\ref{fx1})-(\ref{fx3}) and
$\{g_i|_{X_F}\}_{i=1}^\infty$ is a tight pre-$F$-frame for $X_F$
with respect to $\Theta_F$.
\end{prop}
\bp First observe that for given number sequence $\seqgr[a]$ with
$a_i\geq 1$, $i\newin \mn$, the space $X$ defined by
$$
X:=\left\{f\in X_0 \ : \ a_i \<f,e_i\>_0 \in \ell^2, \ \
\<f,h\>_X:=\sumii a_i^2 \<f,e_i\>_0 \,\<e_i, h\>_0\right\},
$$
is a Hilbert space satisfying $\|\cdot\|_{X_0}\leq \|\cdot\|_X$
and having $\{\frac{1}{a_i}\,e_i\}_{i=1}^\infty$ as an orthonormal
basis. For every $s\newin\mn$, apply this to the sequence
$\{a_{i,s}\}_{i=1}^\infty$ and to the space $X_s$.

Fix $s\newin\mn$. The inequalities  $a_{i,s}\leq a_{i,s+1}$,
$i\newin\mn$, imply that $X_{s+1}\subseteq X_s$ and thus
(\ref{fx1}) holds. Denote $z_{i,s}:=e_i/a_{i,s}$, $i\newin\mn$. By
definition, we have
\begin{equation}\label{econd}
 \<f,z_{i}\>_s= a_{i,s}\, \<f, e_i\>_0, \ f\in X_s,\ i\in\mn.
  \end{equation}
Thus, for $f\in X_{s+1}$ one has
$$\|f\|_s^2=\sumii |\<f,z_{i,s} \>_s |^2=\sumii a_{i,s}^2 |\<f,e_{i} \>_0 |^2\leq \sumii a_{i,s+1}^2 |\<f,e_{i} \>_0 |^2 =
\|f\|_{s+1}^2$$ and hence (\ref{fx2}) holds.

It is clear that $e_i\newin X_F$, $i\newin\mn$. Moreover, the
linear span of $\seqgr[e]$ (which is a subset of $X_F$) is dense
in $X_s$. Therefore, (\ref{fx3}) also holds.

Let us now show that $(\mathcal{A}_1^\prime)$ is fulfilled. Take
$c\newin\Theta_s$, $d\newin\Theta_s$, $f\newin M_s^c$, $h\newin
M_s^d$. We know that $|c_i|\leq |g_i(f)|$, $i\newin \mn$, and
$$ \sum_{i=1}^\infty a_{i,s}^2|g_{i+1}(f)|^2
=\sum_{i=1}^\infty a_{i,s}^2|\<f,e_i\>_0|^2=\sum_{i=1}^\infty
|\<f,z_i\>_s|^2= ||f||_s^2<\infty,
$$ which implies that
\begin{equation} \label{shod} \sum_{i=1}^\infty a_{i,s}^2|c_{i+1}|^2 < \infty;
\ \mbox{similarly,}\ \sum_{i=1}^\infty a_{i,s}^2|d_{i+1}|^2 <
\infty.
\end{equation}
Let $$m:= \max(|c_1+d_1|, |c_2+d_2|), \ \widetilde{c}:=\max(|c_1|,
|c_2|), \ \widetilde{d}:=\max(|d_1|, |d_2|).$$
 It follows from (\ref{shod}) that
\begin{equation}\label{rfh}
r^{f,h}:=ma_{1,s}z_1+|c_3 + d_3|a_{2,s}z_2 + |c_4+d_4|
a_{3,s}z_{3}+\ldots
\end{equation}
 is an element of $X_s$.
By (\ref{econd}) and (\ref{rfh}), it follows that
\begin{eqnarray*}
g_i(r^{f,h})=\<r^{f,h},e_{i-1}\>_0= \frac{1}{a_{i-1,s}}
\<r^{f,h},z_{i-1}\>_s =
|c_i+d_i|, \; i\geq 3,\\
 g_i(r^{f,h})=\<r^{f,h},e_1\>_0 = \frac{1}{a_{1,s}} \<r^{f,h},z_{1}\>_s=m\geq |c_i+d_i|,\; i=1,2.
\end{eqnarray*}
This implies that $r^{f,h}\in M_s^{c+d}$. In a similar way as in
Proposition \ref{primer1} we obtain that
\begin{eqnarray*}
\|r^{f,h}\|_s &=&
\left\|\left\{a_{1,s}m, a_{2,s}|c_3+d_3|, a_{3,s}|c_4+d_3|, \ldots  \right\}\right \|_{\ell^2}\\
&\leq & \left\|\left\{a_{1,s}\widetilde{c}+a_{1,s}\widetilde{d},
a_{2,s}|c_3|+a_{2,s}|d_3|, a_{3,s}|c_4|+a_{3,s}|d_4|,
\ldots \right\} \right\|_{\ell^2}\\
& \leq & \left\|\left\{a_{1,s}\widetilde{c}, a_{2,s}|c_3|,
a_{3,s}|c_4|, \ldots
\right\}\right\|_{\ell^2} \\
 & &+ \left\|\left\{a_{1,s}\widetilde{d},a_{2,s}|d_3|,a_{3,s}|d_4|, \ldots \right\}
\right\|_{\ell^2}\\
& \leq & \left\|\left\{a_{i,s}\<f,e_i\>_0\right\}_{i=1}^\infty\right\|_{\ell^2}  +\left\|\left\{a_{i,s}\<h,e_i\>_0 \right\}_{i=1}^\infty \right\|_{\ell^2}\\
& = &
\left\|\left\{\<f,z_i\>_s\right\}_{i=1}^\infty\right\|_{\ell^2}
+\left\|\left\{\<h,z_i\>_s \right\}_{i=1}^\infty
\right\|_{\ell^2}=\|f\|_s+\|h\|_s.
\end{eqnarray*}

Let us prove that $(\mathcal{A}_2^\prime)$ holds. Take
$c=\{c_1,c_2,c_3,\ldots \}\in \Theta_s$. Fix $\varepsilon>0$ and
by (\ref{shod}), find $k\in\mn$, $k>1,$ such that
$\sum_{i=k}^\infty a_{i,s}^2|c_{i+1}|^2<\varepsilon$. Let
$$b:=\{\underbrace{0, \ldots, 0}_{k-1}, a_{k,s}\, |c_{k+1}|, a_{k+1,s}\,|c_{k+2}|, a_{k+2,s}\,|c_{k+3}|,\ldots \}.$$
Since $b\newin\ell^2$, there exists $h\newin X_s$ such that
$b_i=\<h, z_i\>_s$, $i\newin\mn$, and
 $$\|h\|_s^2=\sumii |b_i|^2=
\sum_{i=k}^\infty a_{i,s}^2|c_{i+1}|^2<\varepsilon.$$ Moreover,
$$|c^{(k)}_1|=0\leq |g_{1}(h)| \ \, \mbox{and} \ \,
|c^{(k)}_i|=\frac{b_{i-1}}{a_{i-1,s}}=\frac{\<h,z_{i-1}\>_s}{a_{i-1,s}}=
g_{i}(h), \ \ i\geq 2,$$ which implies that $h\in M_s^{c^{(k)}}$.

Let us now prove that $(\mathcal{A}_3^\prime)$ is fulfilled.
Consider $f\newin X_s$ and take arbitrary $\widetilde{f}\in
M_s^{\{g_i(f)\}_{i=1}^\infty}$, i.e. $|g_i(f)|\leq
|g_i(\widetilde{f})|, i\in\mn$. Using (\ref{econd}), we obtain
$$
\|f\|_s^2=\sumii a_{i,s}^2 \,|\<f, e_i\>_0|^2 = \sumii a_{i,s}^2\,
|g_{i+1}(f)|^2 \leq  \sumii  a_{i,s}^2\,
|g_{i+1}(\widetilde{f})|^2=\|\widetilde{f}\|_s^2.
$$
Therefore $(\mathcal{A}_3^\prime)$ holds with $A_s=1$.

Now Corollary \ref{c1} implies that $\{g_i|_{X_F}\}_{i=1}^\infty$
is a tight pre-$F$-frame for $X_F$ with respect to $\Theta_F$. \ep

\v We finish the paper with the following open problem:

{\bf Problem.} Find conditions on $\{X_s\}_{s\in\mn_0}$ and
$\seqgr[g]$, implying existence of an $F$-bounded projection from
$\Theta_F$ onto  $\{\{g_i(f)\}_{i=1}^\infty \, : \, f\in X_F\}$,
which would imply that $\seqgr[g]$ is an $F$-frame for $X_F$ with
respect to $\Theta_F$.

\noindent
S. Pilipovi\'c\\
Department of Mathematics and Informatics\\ University of Novi Sad\\ Trg D.\ Obradovi\'ca 4\\
 21000 Novi Sad, Serbia\\
pilipovic@im.ns.ac.yu

\vspace{.1in} \noindent
D.\,T. Stoeva\\
Department of Mathematics\\ University of Architecture, Civil Engineering and Geodesy\\
 Blvd
Christo Smirnenski 1\\ 1046 Sofia, Bulgaria\\
stoeva\_\,fte@uacg.bg

\end{document}